\documentclass{article}

\usepackage{etex}
\usepackage{graphicx}
\usepackage{latexsym}
\usepackage{amsfonts}
\usepackage{amssymb}
\usepackage{amsmath}
\usepackage{verbatim}
\usepackage{url}
\usepackage[standard, thmmarks]{ntheorem}
\usepackage[margin=3cm]{geometry}
\usepackage{fancyhdr}

\newtheorem{ax}{Axiom}
\newtheorem*{ax0}{Axiom}
\newtheorem{thm}{Theorem}[section]
\newtheorem{cor}[thm]{Corollary}
\newtheorem{lem}[thm]{Lemma}
\newtheorem{prop}[thm]{Proposition}

\newtheorem{defn}[thm]{Definition}

\newcommand{\R}{\mathbb R}

\newcommand{\Z}{\mathbb Z}

\newcommand{\Q}{\mathbb Q}

\newcommand{\F}{\mathcal F}

\newcommand{\To}{\longrightarrow}

\input xy
\xyoption{all}

\begin{document}

\title{Sutured TQFT, torsion, and tori}

\author{Daniel V. Mathews}%

\date{}

\maketitle

\begin{abstract}
We use the theory of sutured TQFT to classify contact elements in the sutured Floer homology, with $\Z$ coefficients, of certain sutured manifolds of the form $(\Sigma \times S^1, F \times S^1)$ where $\Sigma$ is an annulus or punctured torus. Using this classification, we give a new proof that the contact invariant in sutured Floer homology with $\Z$ coefficients of a contact structure with Giroux torsion vanishes. We also give a new proof of Massot's theorem that the contact invariant vanishes for a contact structure on $(\Sigma \times S^1, F \times S^1)$ described by an isolating dividing set.

\end{abstract}

\tableofcontents

\section{Introduction}

\subsection{Torsion and isolation}

In this paper we give a new proof of an important fact about contact topology and Heegaard Floer homology.

In 3-dimensional contact topology, the idea of \emph{torsion} first arose in the work of Giroux \cite{GiPlusMoins, Gi99, Gi00}. A contact manifold has $(2\pi)$-\emph{torsion} if it admits a contact embedding of $(T^2 \times [0,1], \xi)$ where $\xi = \ker (\cos (2\pi t) \; dx - \sin (2\pi t) \; dy )$, and $((x,y),t)$ are coordinates on $T^2 \times [0,1]$. Torsion has played an important role in classifying contact structures on 3-manifolds \cite{Colin01_torsion, Colin01_structures, Colin_Giroux_Honda03, HKM02, Colin_Giroux_Honda09} and fillability \cite{Gay06}.

A contact structure $\xi$ on a 3-manifold also gives rise to a \emph{contact invariant} or \emph{contact elements} $c(\xi)$ in the Heegaard Floer homology of Ozsv\'{a}th--Szab\'{o} \cite{OS04Closed} or its younger sibling, sutured Floer homology, defined by Juh\'{a}sz \cite{Ju06}. Both in the closed case \cite{OSContact, HKMContClass} and the sutured case \cite{HKM09}, contact elements possess various natural properties. It was conjectured by Ghiggini in \cite{Ghiggini06Infinite} that the contact invariant vanishes for a torsion contact structure. After partial results in \cite{Ghiggini06Fillability, Lisca_Stipsicz07}, Ghiggini--Honda--Van Horn Morris gave a proof in \cite{GHV}, and Massot gave another proof in \cite{Massot09}. In this paper we give a third proof, based on the theory of sutured TQFT, developed in \cite{Me09Paper} and \cite{Me10_Sutured_TQFT}.

\begin{thm}
\label{thm:torsion_vanishing}
Let $\xi$ be a contact structure with $(2\pi)$-torsion on a closed 3-manifold $M$ or a balanced sutured 3-manifold $(M,\Gamma)$, with contact invariant $c(\xi) \subset \widehat{HF}(-M)$ or $SFH(-M,-\Gamma)$ respectively (with $\Z$ coefficients). Then $c(\xi) = \{0\}$.
\end{thm}

In \cite{Massot09}, Massot also proved a conjecture of Honda--Kazez--Mati\'{c} \cite{HKM08}. Take a surface $\Sigma$ with nonempty boundary, and form the balanced sutured 3-manifold $(\Sigma \times S^1, F \times S^1)$ where $F$ is a finite subset of $\partial \Sigma$. Drawing a dividing set $\Gamma$ on $(\Sigma, F)$ gives a contact structure $\xi_\Gamma$ on $(\Sigma \times S^1, F \times S^1)$; and in fact the set of dividing sets (up to isotopy) without homotopically trivial closed curves on $(\Sigma, F)$ is bijective with the set of tight contact structures (up to isotopy) on $(\Sigma \times S^1, F \times S^1)$; $\{\Gamma\} \leftrightarrow \{\xi_\Gamma\}$ is a bijection \cite{GiBundles, Hon00II}. A dividing set $\Gamma$ on $(\Sigma, F)$ is \emph{isolating} if $\Sigma \backslash \Gamma$ has a component not intersecting $\partial \Sigma$. The conjecture of Honda--Kazez--Mati\'{c}, proved by Massot, is as follows. 
\begin{thm}
\label{thm:isolating_vanishing}
Consider contact structures $\xi$ on $(\Sigma \times S^1, F \times S^1)$ described by dividing sets $\Gamma$ on $(\Sigma, F)$. In sutured Floer homology over $\Z$ coefficients, the following are equivalent:
\begin{enumerate}
\item
$c(\xi) \neq 0$.
\item
$c(\xi)$ is primitive.
\item
$\Gamma$ is not isolating.
\end{enumerate}
\end{thm}
In this paper we shall give another proof of this result, again based on sutured TQFT.

On manifolds $(\Sigma \times S^1, F \times S^1)$, a torsion contact structure is described by an isolating dividing set, so for this class of manifolds theorem \ref{thm:isolating_vanishing} is a generalisation of theorem \ref{thm:torsion_vanishing}.  In \cite{HKM08}, Honda--Kazez--Mati\'{c} proved this result over $\Z_2$ coefficients and proved $(iii) \Rightarrow (ii) \Rightarrow (i)$ over $\Z$ coefficients. In \cite{Massot09}, Massot proved $(i) \Rightarrow (iii)$ over $\Z$ coefficients, completing the proof.

\subsection{Sutured TQFT marches on}

As mentioned above, the results in this paper are obtained through sutured TQFT, as developed in \cite{Me09Paper, Me10_Sutured_TQFT}, which arises from from the TQFT properties of sutured Floer homology introduced in \cite{HKM08}. The case of sutured manifolds of the form $(\Sigma \times S^1, F \times S^1)$ can be regarded as a dimensionally-reduced case, with some properties analogous to a (1+1)-dimensional TQFT. (Although, as mentioned in \cite{Me10_Sutured_TQFT}, in some ways it is better considered a ``(2+1=2)-dimensional TQFT''.) 

The theory of sutured TQFT abstracts from this dimensionally reduced case of sutured Floer homology, to give an axiomatically defined theory formally free of holomorphic curves and contact geometry. In essence, sutured TQFT associates to a surface with some decorations on the boundary (a \emph{sutured background surface}) $(\Sigma, F)$  an abelian group $V(\Sigma, F)$; and to certain sets of curves $\Gamma$ drawn on the surface (\emph{sutures}), associates a subset of \emph{suture elements} $c(\Gamma) \subset V(\Sigma, F)$. These associations are required to be natural with respect to gluing operations, and various other axioms, inspired by TQFT, sutured Floer homology and contact geometry, are imposed in order to give structure to the theory. We shall state the axioms below in section \ref{sec:sutured_TQFT}.

In \cite{Me10_Sutured_TQFT} we showed that the sutured Floer homology of manifolds $(\Sigma \times S^1, F \times S^1)$ forms a sutured TQFT, essentially taking $V(\Sigma,F)$ to be $SFH(\Sigma \times S^1, F \times S^1)$ and $c(\Gamma)$ to be contact elements. Thus, any result about sutured TQFT also gives a result about sutured Floer homology. Since results can be obtained in sutured TQFT by purely combinatorial or topological methods, without use of holomorphic curves or even contact geometry, we can obtain results about sutured Floer homology with proofs that are ``holomorphic curve free'' or ``contact geometry free''.

In \cite{Me09Paper} we computed sutured TQFT of discs $D^2$ over $\Z_2$ in detail, obtaining a rich algebraic and combinatorial structure; and thus for contact elements in the sutured Floer homology with $\Z_2$ coefficients of solid tori with longitudinal sutures. In \cite{Me10_Sutured_TQFT} we extended these results to $\Z$ coefficients, obtaining an even richer and more general algebraic and combinatorial structure, which can be described as a Fock space of two non-commuting particles.

In this paper, we perform further computations in sutured TQFT. Having completely determined the structure of sutured TQFT of discs in \cite{Me10_Sutured_TQFT}, we move on to annuli and tori. We compute the structure of sutured TQFT of annuli and once-punctured tori with certain simple boundary markings. We give a complete classification of suture elements, for these background surfaces. This then immediately gives a classification of contact elements in $SFH$ of the corresponding manifolds. The structure of contact elements in these cases is quite interesting; we are able to take coordinates on $SFH$ in such a way that the slope of the dividing set gives the coordinates of the corresponding contact element. The classification over annuli is theorem \ref{thm:annulus_contact_elt_classification}; it was also effectively obtained by Massot in \cite{Massot09}.
\begin{thm}
\label{thm:annulus_contact_elt_classification}
Let $\Sigma$ be an annulus and $F$ two points on each boundary component of $\Sigma$. Then $SFH(-\Sigma \times S^1, -F \times S^1) \cong \Z^4$, splitting as summands $\Z \oplus \Z^2 \oplus \Z$, corresponding to contact elements with euler class $-2,0,2$. The nonzero contact elements, and corresponding dividing sets $\Gamma$ on $(\Sigma, F)$, are precisely as follows.
\begin{enumerate}
\item
The only contact structure with euler class $-2$ and nonzero contact element corresponds to $\Gamma$ consisting of two boundary-parallel arcs on $\Sigma$, each enclosing a negative disc; the contact element is $\pm 1 \in \Z$.
\item
The only contact structures with euler class $0$ and nonzero contact element correspond to dividing sets $\Gamma$ as follows.
\begin{enumerate}
\item
A closed loop around the core of the annulus, and two boundary-parallel arcs; the contact element is $\pm(0,1) \in \Z^2$.
\item
Two parallel arcs between boundary components, traversing the core of $\Sigma$ $n$ times, i.e. having slope $n/1$; the contact element is $\pm(1,n) \in \Z^2$.
\end{enumerate}
\item
The only contact structure with euler class $2$ and nonzero contact element corresponds to $\Gamma$ consisting of two boundary-parallel arcs on $\Sigma$, each enclosing a positive disc; the contact element is $\pm 1 \in \Z$.
\end{enumerate}
\end{thm}

Over the punctured torus, we obtain an even more striking result relating slopes of dividing sets and coordinates of contact elements, closely related to the Farey graph, as we discuss in section \ref{sec:classification}.
\begin{thm}
\label{thm:punctured_torus_contact_elt_classification}
Let $\Sigma$ be a once punctured torus and $F$ two points on its boundary. Then $SFH(-\Sigma \times S^1, -F \times S^1) \cong \Z^4$, splitting as summands $\Z \oplus \Z^2 \oplus \Z$, corresponding to contact elements with euler class $-2,0,2$. The nonzero contact elements, and corresponding dividing sets $\Gamma$ on $(\Sigma, F)$, are precisely as follows.
\begin{enumerate}
\item
The only contact structure with euler class $-2$ and nonzero contact element corresponds to $\Gamma$ consisting of one boundary-parallel arc enclosing a negative disc; the contact element is $\pm 1 \in \Z$.
\item
The only contact structures with euler class $0$ and nonzero contact element are contact structures $\xi_{q/p}$ corresponding to dividing sets $\Gamma_{q/p}$, where $\Gamma_{q/p}$ consists of an arc and a closed loop of slope $q/p$. We may choose bases for $H_1(T_1)$ (to define slope) and $SFH(\Sigma \times S^1, F \times S^1)$ so that $c(\xi_{q/p}) = \pm(p,q)$.
\item
The only contact structure with euler class $2$ and nonzero contact element corresponds to $\Gamma$ consisting of one boundary-parallel arc enclosing a positive disc; the contact element is $\pm 1 \in \Z$.
\end{enumerate}
\end{thm}
So far as we know, theorem \ref{thm:punctured_torus_contact_elt_classification} is new.

On the way to these classification results, we show that in sutured TQFT, sets of sutures with \emph{torsion}, corresponding to torsion contact structures, have zero suture element. More generally, we show that \emph{isolating} sets of sutures, as defined above, also have zero suture elements (sutures with torsion are isolating).
\begin{thm}
\label{thm:STQFT_isolating_vanishing}
In any sutured TQFT, if $\Gamma$ is an isolating set of sutures on a sutured background $(\Sigma, F)$, then $c(\Gamma) = \{0\}$.
\end{thm}
From this theorem, a proof of theorem \ref{thm:isolating_vanishing} is essentially immediate.

\subsection{What this paper does}

As the discussion above indicates, this paper is primarily concerned with computing suture elements in sutured TQFT. We classify suture elements in $V(\Sigma, F)$ where $\Sigma$ is an annulus or once punctured torus, and $F$ is minimal as in theorems \ref{thm:annulus_contact_elt_classification} and \ref{thm:punctured_torus_contact_elt_classification}.

In section \ref{sec:story_so_far} we recall sutured TQFT and previous results on discs. In section \ref{sec:preliminaries} we make some preliminary observations, including various suture elements of chord diagrams based on previous work. In section \ref{sec:annuli} we consider annuli, proving that torsion sutures give zero suture element and classifying suture elements. In section \ref{sec:punctured_tori} we consider once-punctured tori and classify suture elements; along the way we show isolating sutures give zero suture element. Finally, in section \ref{sec:SFH} we consider contact elements in sutured Floer homology. Most of the main theorems follow immediately. However the proof of theorem \ref{thm:torsion_vanishing}, being a result in contact topology, requires us to make some extra considerations.

This paper relies heavily on previous work in \cite{Me10_Sutured_TQFT} and \cite{Me09Paper}. Results about the sutured TQFT of surfaces are obtained from results about the sutured TQFT of discs, gluing up discs into more complicated surfaces. We attempt to gather what we need in sections \ref{sec:story_so_far}--\ref{sec:preliminaries}, and give references to those papers where we can, but some level of familiarity with them must be assumed.

\subsection{Acknowledgments}

This paper was partially written during the author's visit to the Mathematical Sciences Research Institute in March 2010, and during the author's postdoctoral fellowship at the Universit\'{e} de Nantes, supported by the ANR grant ``Floer power''.

\section{The story so far}
\label{sec:story_so_far}

\subsection{Sutured TQFT}
\label{sec:sutured_TQFT}

In \cite{Me10_Sutured_TQFT} we defined \emph{sutured topological quantum field theory}. This theory is designed to be an axiomatic version of the (1+1)-dimensional dimensionally-reduced TQFT-like structure defined by Honda--Kazez--Mati\'{c} in \cite{HKM08}, and to describe contact elements in the sutured Floer homology of manifolds of the form $(\Sigma \times S^1, F \times S^1)$.

We recall some definitions: see \cite[section 3.1]{Me10_Sutured_TQFT} for details. A \emph{sutured surface} $(\Sigma, \Gamma)$ is a compact oriented surface $\Sigma$, possibly disconnected, each component with nonempty boundary, with $\Gamma \subset \Sigma$ a properly embedded oriented 1-submanifold satisfying the following property: $\Sigma \backslash \Gamma = R_+ \cup R_-$, where the $R_\pm$ are oriented as $\pm \Sigma$, and where $\partial R_\pm \cap \Gamma = \pm \Gamma$ as oriented 1-manifolds.

A \emph{sutured background surface} (or simply \emph{sutured background}) $(\Sigma, F)$ is a compact oriented surface $\Sigma$, possibly disconnected, each component with nonempty boundary, together with a finite set of signed points $F \subset \partial \Sigma$, such that $\partial \Sigma \backslash F = C_+ \cup C_-$, where $C_\pm$ are arcs oriented as $\pm \partial \Sigma$, and $\partial C_+ = \partial C_- = F$ as sets of signed points. (Hence each boundary component $C$ of $\Sigma$ has a positive even number of points of $F$, which cut it alternately into arcs of $C_+$ and $C_-$.)

A \emph{set of sutures} $\Gamma$ on a sutured background $(\Sigma,F)$ is an an oriented properly embedded 1-submanifold of $\Sigma$ such that $\partial \Gamma = \partial \Sigma \cap \Gamma = F$ and such that $(\Sigma, \Gamma)$ is a sutured surface, with $\partial R_\pm = \pm \Gamma \cup C_\pm \cup F$. A set of sutures $\Gamma$ has an \emph{Euler class} defined by $e(\Gamma) = \chi(R_+) - \chi(R_-)$.

A sutured background $(\Sigma, F)$ may be \emph{glued}. Let $\partial \Sigma \backslash F = C_+ \cup C_-$ as above. Consider two disjoint 1-manifolds $G_0, G_1 \subseteq \partial \Sigma$, and a homeomorphism $\tau: G_0 \stackrel{\cong}{\To} G_1$ which identifies marked points and positive/negative arcs, $G_0 \cap F \stackrel{\cong}{\to} G_1 \cap F$, $G_0 \cap C_\pm \stackrel{\cong}{\to} G_1 \cap C_\pm$. Then we may glue $(\Sigma, F)$ along $\tau$ and obtain a surface $\#_\tau (\Sigma, F)$. If there remain marked points on each boundary component then $\#_\tau(\Sigma,F)$ is also a sutured background surface and we call $\tau$ a \emph{sutured gluing map}. If $\Gamma$ is a set of sutures on $(\Sigma,F)$ then a sutured gluing map gives a glued set of sutures $\#_\tau \Gamma$ on $\#_\tau (\Sigma,F)$.

For a disc $D^2$, a sutured background is specified simply by giving the number of points in $F$. We define $(D^2, F_n)$ to be the sutured background disc with $|F_n|=2n$. The simplest case $(D^2,F_1)$ is called the \emph{vacuum background}; it only has one set of sutures without closed component, obtained by joining the two points of $F_1$; this set of sutures is called the \emph{vacuum} $\Gamma_\emptyset$.

Sutured TQFT is defined by the following set of axioms. In \cite{Me10_Sutured_TQFT} we introduced these axioms carefully, giving a rationale for each, and explaining possible variations in detail.
\begin{ax}
\label{ax:1}
To each sutured background surface $(\Sigma, F)$, assign an abelian group $V(\Sigma, F)$, depending only on the homeomorphism type of the pair $(\Sigma,F)$.
\end{ax}

\begin{ax}
\label{ax:2}
To a set of sutures $\Gamma$ on $(\Sigma,F)$, assign a subset of \emph{suture elements} $c(\Gamma) \subset V(\Sigma, F)$, depending only on the isotopy class of $\Gamma$ relative to boundary.
\end{ax}

\begin{ax}
\label{ax:3}
For a sutured gluing map $\tau$ of a sutured background surface $(\Sigma, F)$, assign a collection of linear maps $\Phi_\tau^i \; : \; V(\Sigma, F) \To V(\#_\tau(\Sigma, F))$.
\end{ax}

Axiom 3 also has an equivalent formulation 3' in terms of \emph{inclusions}: 
\begin{ax0}[3']
To an inclusion $(\Sigma_{in}, F_{in}) \stackrel{\iota}{\hookrightarrow} (\Sigma_{out},F_{out})$ of sutured background surfaces, with $\Sigma_{in}$ lying in the interior of $\Sigma_{out}$, together with $\Gamma$ a set of sutures on $( \Sigma_{out} \backslash \Sigma_{in}, F_{in} \cup F_{out})$, assign a collection of linear maps $\Phi_{\iota,\Gamma}^i \; : \; V(\Sigma_{in}, F_{in}) \To V(\Sigma_{out}, F_{out})$.
\end{ax0}

\begin{ax}
\label{ax:4}
For a finite disjoint union of sutured background surfaces $\sqcup_i (\Sigma_i, F_i)$, 
\[
V(\sqcup_i (\Sigma_i, F_i)) = \otimes_i V(\Sigma_i, F_i).
\]
\end{ax}

\begin{ax}
\label{ax:5}
If $\Gamma$ is a set of sutures on $(\Sigma,F)$ and $\tau$ is a gluing of $(\Sigma,F)$ then each $\Phi_\tau^i$ takes suture elements to suture elements surjectively, $c(\Gamma) \to c(\#_\tau \Gamma)$.
\end{ax}

Corresponding to axiom 3' in terms of inclusions is axiom 5'. Axioms 3 and 5 are equivalent to axioms 3' and 5'.
\begin{ax0}[5']
If $\Gamma_{in}$ is a set of sutures on $(\Sigma_{in}, F_{in})$, let $\Gamma_{out} = \Gamma_{in} \cup \Gamma$ be the corresponding set of sutures on $(\Sigma_{out}, F_{out})$. Then each $\Phi_{\iota,\Gamma}^i$ maps takes suture elements to suture elements surjectively, $c(\Gamma_{in}) \to c(\Gamma_{out})$.
\end{ax0}

\begin{ax}
\label{ax:6}
If $\Gamma$ contains a closed contractible loop then $c(\Gamma) = \{0\}$.
\end{ax}

\begin{ax}
\label{ax:7}
$V(D^2, F_1) = \Z$ and $c(\Gamma_\emptyset) \subseteq \{-1,1\}$.
\end{ax}

\begin{ax}
\label{ax:8}
Every $V(\Sigma,F)$ is spanned by suture elements.
\end{ax}

\begin{ax}
\label{ax:nondegeneracy}
Suppose two elements $x,y \in V(D^2,F_n)$ have the following property: for any set of sutures $\Gamma$ on $(D^2,F_n)$, there exists $c \in c(\Gamma)$ such that $\langle x | c \rangle = \pm \langle y | c \rangle$. Then $x = \pm y$.
\end{ax}

In \cite{Me10_Sutured_TQFT} we also introduced an additional tenth axiom, which is stronger and implies axiom \ref{ax:nondegeneracy}, and almost implies axiom \ref{ax:8} also.
\begin{ax}
\label{ax:10}
Let $\tau$ be a sutured gluing map on $(\Sigma,F)$, identifying two disjoint arcs $\gamma, \gamma'$ on $\partial \Sigma$. Suppose that $|\gamma \cap F| = |\gamma' \cap F| = 1$. Then any gluing map $\Phi_\tau$ associated to $\tau$ is an isomorphism.
\end{ax}

We shall assume all 10 axioms for the purposes of this paper.

\subsection{Results for discs}

In \cite{Me09Paper, Me10_Sutured_TQFT} we investigated the structure of sutured TQFT, in particular of discs. We found that this was isomorphic to a ``Fock space of two non-commuting particles'' $\F$. We summarise some of this structure as follows, but refer there for details.

In general, given a set of sutures $\Gamma$ on a sutured background $(\Sigma, F)$, $c(\Gamma)$ is of the form $\{\pm x\}$, so unless $x$ is torsion or $0$, $c(\Gamma)$ has cardinality two and $c(\Gamma)$ has a sign ambiguity.

A set of sutures $\Gamma$ on a sutured background disc $(D^2, F_n)$ always either has a closed contractible component, so that $c(\Gamma) = \{0\}$, or is a \emph{chord diagram}, a properly embedded collection of disjoint arcs in the disc (up to homotopy). Rotating a chord diagram in general gives a distinct chord diagram; in order to keep track of points, we label the $2n$ marked points on the boundary. We choose a basepoint, which is labelled $0$, and then then $2n$ points are labelled clockwise mod $2n$.

We can define $\F$ as $\Z[x,y]$, where $x$ and $y$ are non-commuting variables; so elements of $\F$ are $\Z$-linear combinations of words in $\{x,y\}$. This forms a (non-commutative) ring, which is bi-graded by degree $n_x,n_y$ in $x,y$ respectively, $\F = \oplus_{n_x, n_y} \F_{n_x,n_y}$; alternatively, letting $n=n_x + n_y$ and $e= n_y - n_x$, $\F$ we have another bi-grading, and $\F$ decomposes as $\F = \oplus_{n,e} \F_n^e$; here $\F_{n_x,n_y}$ consists of $\Z$-linear combinations of words with $n_x$ $x$'s and $n_y$ $y$'s; and $\F_n^e$ consists of $\Z$-linear combinations of words with $n$ letters and ``charge'' $e$, where an $x$ has charge $-1$ and $y$ has charge $1$. We also define $\F_n = \oplus_e \F_n^e$.

On $\F$ we define a great deal of structure. We define creation and annihilation operators obeying a bi-simplicial structure. On words there is a partial order $\leq$ given by $w_0 \leq w_1$ if $w_1$ can be obtained from $w_0$ by moving $x$'s only to the right (and moving $y$'s to the left). From this we define a bilinear form $\langle \cdot | \cdot \rangle: \F \otimes \F \To \Z$ by $\langle w_0 | w_1 \rangle = 1$ if $w_0 \leq w_1$, and $0$ otherwise. Since $\langle \cdot | \cdot \rangle$ is nondegenerate we can define a ``duality'' operator $H$ defined by $\langle u | v \rangle = \langle v | Hu \rangle$.

We find that this all corresponds precisely to structure in sutured TQFT. Each $V(D^2, F_{n+1}) \cong \F_n$, and so if we set $V(D^2) = \oplus_n V(D^2, F_n)$ then $V(D^2) \cong \F$. Corresponding to the decomposition $\F_n = \oplus_e \F_n^e$, we have the decomposition $V(D^2, F_{n+1}) = \oplus_e V(D^2, F_{n+1})^e$ where the summand $\F_n^e$ corresponds to the subgroup $V(D^2,F_{n+1})^e$ of $V(D^2, F_{n+1})$ spanned by suture elements of sutures with Euler class $e$. We shall identify $V(D^2) = \F$ throughout, writing suture elements in $V(D^2)$ as linear combinations of words in $x$ and $y$.

The creation and annihilation operators on $\F$ correspond to inserting or closing off outermost arcs to sets of sutures on $(D^2, F_n)$. A certain subset of chord diagrams, called \emph{basis} chord diagrams, constructed in \cite{Me09Paper}, correspond to the basis of $\F$, i.e. the words on $\{x,y\}$. The bilinear form $\langle \cdot | \cdot \rangle$ corresponds to ``stacking'' contact structures as defined in \cite{Me09Paper}. The partial order on words tells us the stackability of basis chord diagrams.

There is an operation on sets of sutures $\Gamma$ called (upwards or downwards) \emph{bypass surgery}, obtained by taking a disc $D' \subset \Sigma$ which intersects $\Gamma$ in three parallel arcs, and altering the sutures so as to perform a $60^\circ$ rotation (clockwise or anticlockwise) on $D'$; this corresponds to bypass attachment above or below a convex surface in contact geometry \cite{Gi91, Hon00I}. Sets of sutures related by bypass surgery naturally come in triples (\emph{bypass triples}), and their suture elements sum to zero, when signs are appropriately chosen. Using this property, in \cite{Me09Paper} we showed how a chord diagram on $(D^2, F_n)$ can be decomposed as a linear combination of basis diagrams; in \cite{Me10_Sutured_TQFT} we gave a technique to resolve the sign issues involved.

Thus, a general chord diagram with $n+1$ chords and Euler class $e$ has suture element in $\F_n^e$, which is a linear combination of words in $\{x,y\}$. We showed in \cite{Me10_Sutured_TQFT} that all coefficients arising are $\pm 1$, and for a non-basis chord diagram the coefficients sum to $0$; moreover suture elements of chord diagrams are closed under multiplication. Among the words occurring in the basis decomposition of a suture element of a chord diagram $\Gamma$, there is a first and last word $w_-, w_+$ with respect to $\leq$; any word $w$ occurring in this linear combination satisfies $w_- \leq w \leq w_+$. In fact there is a bijection between comparable pairs of words $(w_- \leq w_+)$ and chord diagrams $\Gamma$.

All of the above is described in detail in \cite{Me09Paper} and \cite{Me10_Sutured_TQFT}.

\section{Preliminaries}
\label{sec:preliminaries}

\subsection{Strong sutured TQFT}
\label{sec:zero}

Having considered discs in \cite{Me10_Sutured_TQFT}, we turn to annuli and punctured tori. We shall use the fact that any connected sutured background surface can be formed by taking a sutured background disc, and performing sutured gluings; so we shall rely heavily on our results for discs in \cite{Me10_Sutured_TQFT}.

The following definitions of \emph{torsion} and \emph{isolating} sutures are central to this paper. Clearly a set of sutures with torsion is isolating.
\begin{defn}
A set of sutures $\Gamma$ on a sutured background $(\Sigma, F)$ has \emph{torsion} if there is a component of $\Sigma \backslash \Gamma$ which is an annulus bounded by two closed curves of $\Gamma$.
\end{defn}

\begin{defn}
A set of sutures $\Gamma$ on a sutured background $(\Sigma, F)$ is \emph{isolating} if there is a component of $\Sigma \backslash \Gamma$ which does not intersect $\partial \Sigma$.
\end{defn}

As we saw in \cite{Me10_Sutured_TQFT}, axioms \ref{ax:1}--\ref{ax:nondegeneracy} determine the sutured TQFT of discs completely. This is clearly not true for higher genus surfaces. A gluing $\tau$ can only increase genus: therefore, for instance, we could simply set $V(\Sigma,F) = 0$ whenever $\Sigma$ has genus at least $1$, set all gluing maps to higher genus surfaces to be $0$, and we would obtain a theory consistent with axioms \ref{ax:1}--\ref{ax:nondegeneracy}. However, under the additional assumption of axiom \ref{ax:10}, the theory is not zero at higher genus; we have a stronger version of sutured TQFT. We assume all axioms \ref{ax:1}--\ref{ax:10} throughout this paper.

Under this stronger set of axioms, we easily obtain some basic results about sutured TQFT. These essentially appeared in \cite{HKM08}.
\begin{prop} \cite[lemma 7.2]{HKM08}
\label{prop:no_torsion}
Let $(\Sigma, F)$ be a sutured background surface with $|F|=2n$. As an abelian group, $V(\Sigma, F) \cong \Z^{2^{n - \chi(\Sigma)}}$.
\end{prop}
In particular, no suture element (except for $0$) has torsion (in the sense of group theory!) in this version of sutured TQFT.

\begin{Proof}
Such a $(\Sigma, F)$ can be constructed by gluing together precisely $n- \chi(\Sigma)$ discs $(D^2, F_2)$ by gluing maps of the type considered in axiom 10. By the isomorphism of \cite{Me10_Sutured_TQFT}, $V(D^2, F_2) \cong \F_1 \cong \Z^2$. By axiom 10, these gluing maps are all isomorphisms, so we obtain an isomorphism
\[
\Z^{2^{n - \chi(\Sigma)}} = \bigotimes_{i=1}^{n-\chi(\Sigma)} \Z^2 = \bigotimes_{i=1}^{n-\chi(\Sigma)} V(D^2,F_2) = V \left( \bigsqcup_{i=1}^{n - \chi(\Sigma)} (D^2, F_2) \right) \cong V(\Sigma,F).
\]
\end{Proof}

\begin{prop} \cite[proposition 7.10]{HKM08}
\label{prop:nonzero}
Let $\Gamma$ be a non-isolating set of sutures on $(\Sigma, F)$. Then $c(\Gamma) \neq 0$ and $c(\Gamma)$ is primitive.
\end{prop}

\begin{Proof}
The result is clearly true for discs, from the above.

We show it is possible to find a properly embedded arc $\gamma$ which is not boundary parallel, transverse to $\Gamma$ and such that $|\gamma \cap \Gamma| \leq 1$. To see this, suppose every non-boundary parallel properly embedded arc $\gamma$ intersects $\Gamma$ in at least two points. Denote the components of $\Sigma \backslash \Gamma$ as $R_i$. Traversing $\gamma$, we pass from one region $R_1$ into another region $R_2$; as $\Gamma$ is non-isolating, $R_2$ intersects $\partial \Sigma$. Let $\alpha$ be an arc in $R_2$ connecting $\partial \Sigma$ to $\gamma$, and splitting $\gamma$ into two sub-arcs $\gamma_1, \gamma_2$. Then $\gamma_1 \cup \alpha$ and $\gamma_2 \cup \alpha$ are both properly embedded arcs with fewer intersections with $\Gamma$ than $\gamma$; by minimality then $\gamma_1 \cup \alpha$ and $\gamma_2 \cup \alpha$ are both boundary-parallel. It follows that $\gamma$ is boundary parallel, a contradiction.

Perturb $\gamma$ near $\partial \Sigma$ is necessary so that $|\gamma \cap \Gamma| = 1$, cut along $\gamma$ and repeat. We eventually cut $\Sigma$ along disjoint properly embedded arcs $\gamma_i$, into discs. But the suture elements on discs are nonzero and primitive, and by axiom 10, gluing along the $\gamma_i$ give isomorphisms on sutured TQFT, so that $c(\Gamma)$ is also nonzero and primitive.
\end{Proof}

\subsection{Menagerie}

\begin{figure}
\centering
\includegraphics[scale=0.4]{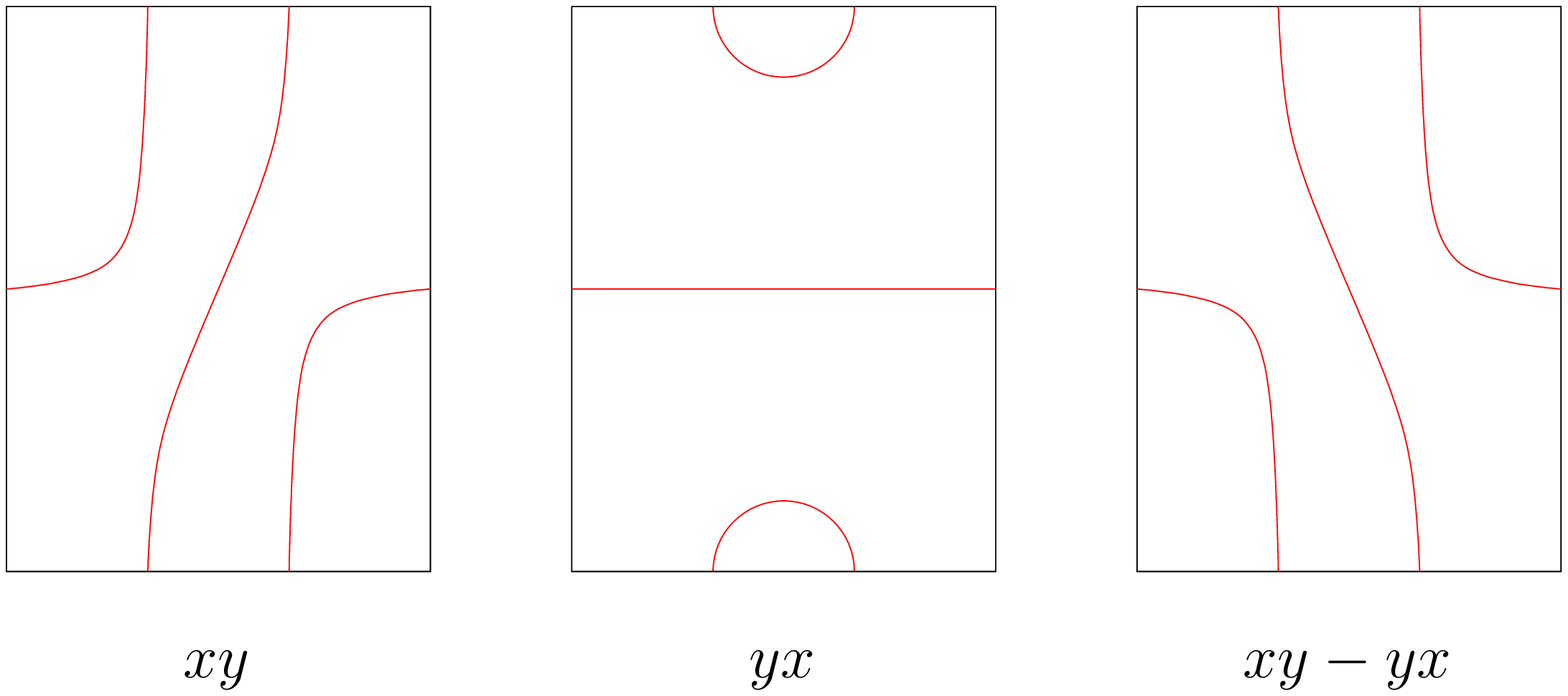}
\caption{Chord diagrams with $3$ chords.} \label{fig:22}
\end{figure}

\begin{figure}
\centering
\includegraphics[scale=0.4]{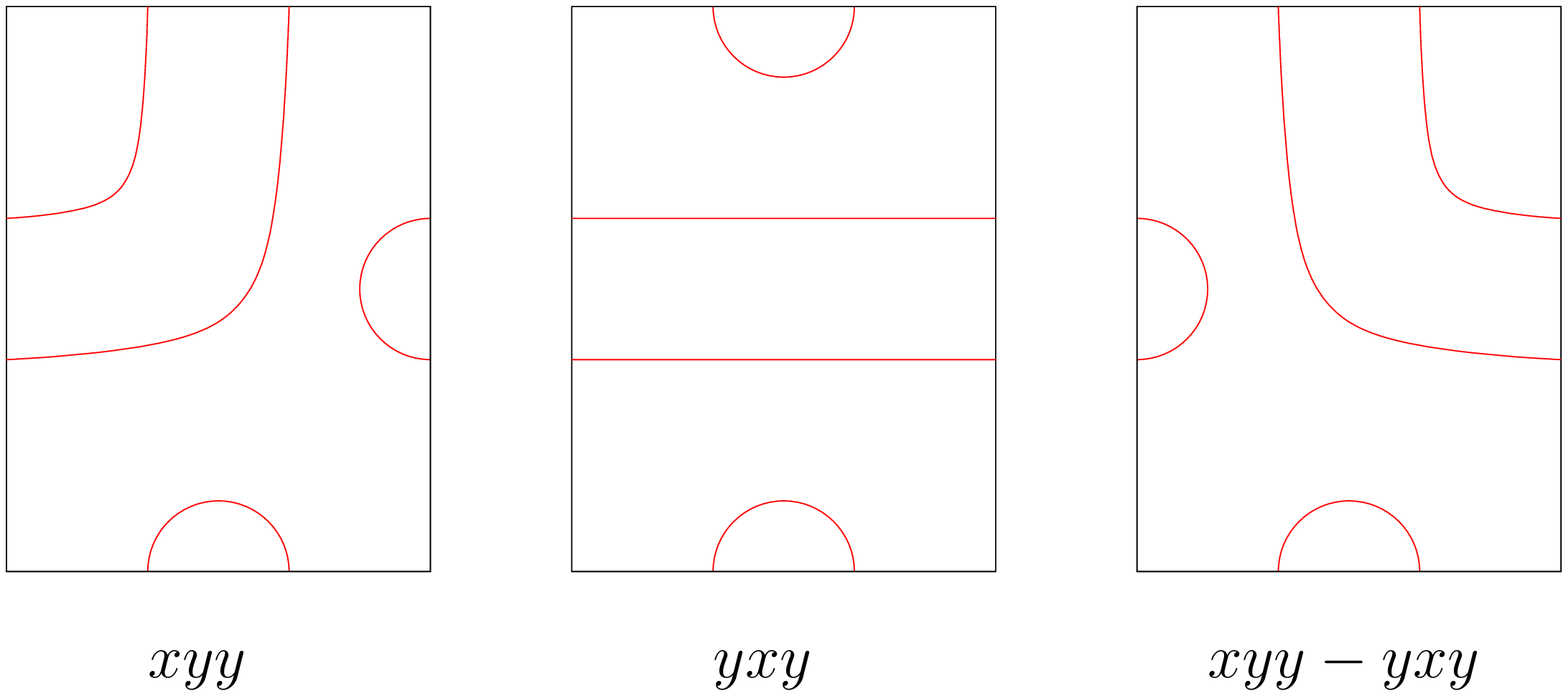}
\caption{Some chord diagrams with $4$ chords.} \label{fig:23}
\end{figure}

\begin{figure}
\centering
\includegraphics[scale=0.4]{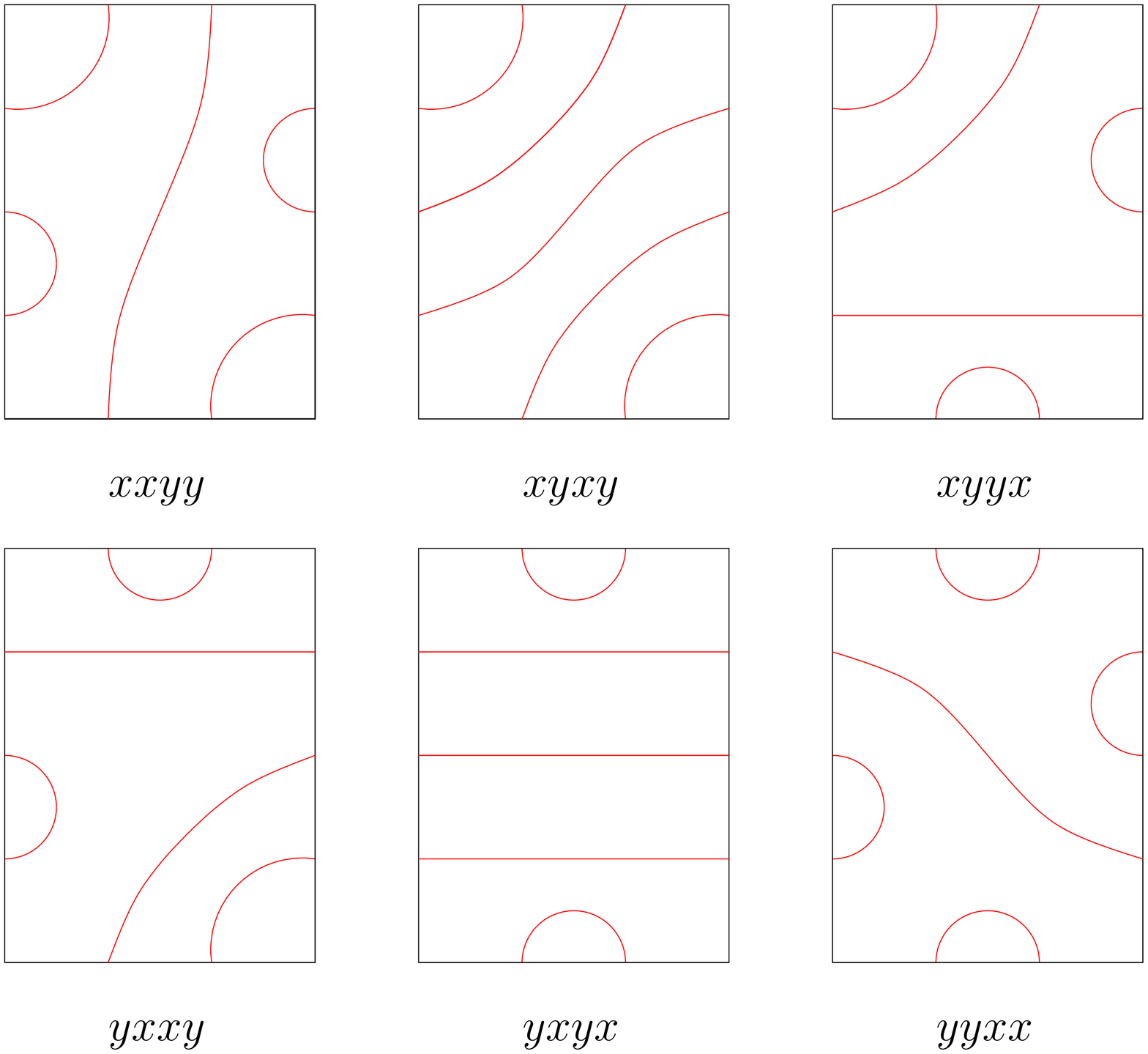}
\caption{Some basis chord diagrams with $5$ chords.} \label{fig:24}
\end{figure}

\begin{figure}
\centering
\includegraphics[scale=0.4]{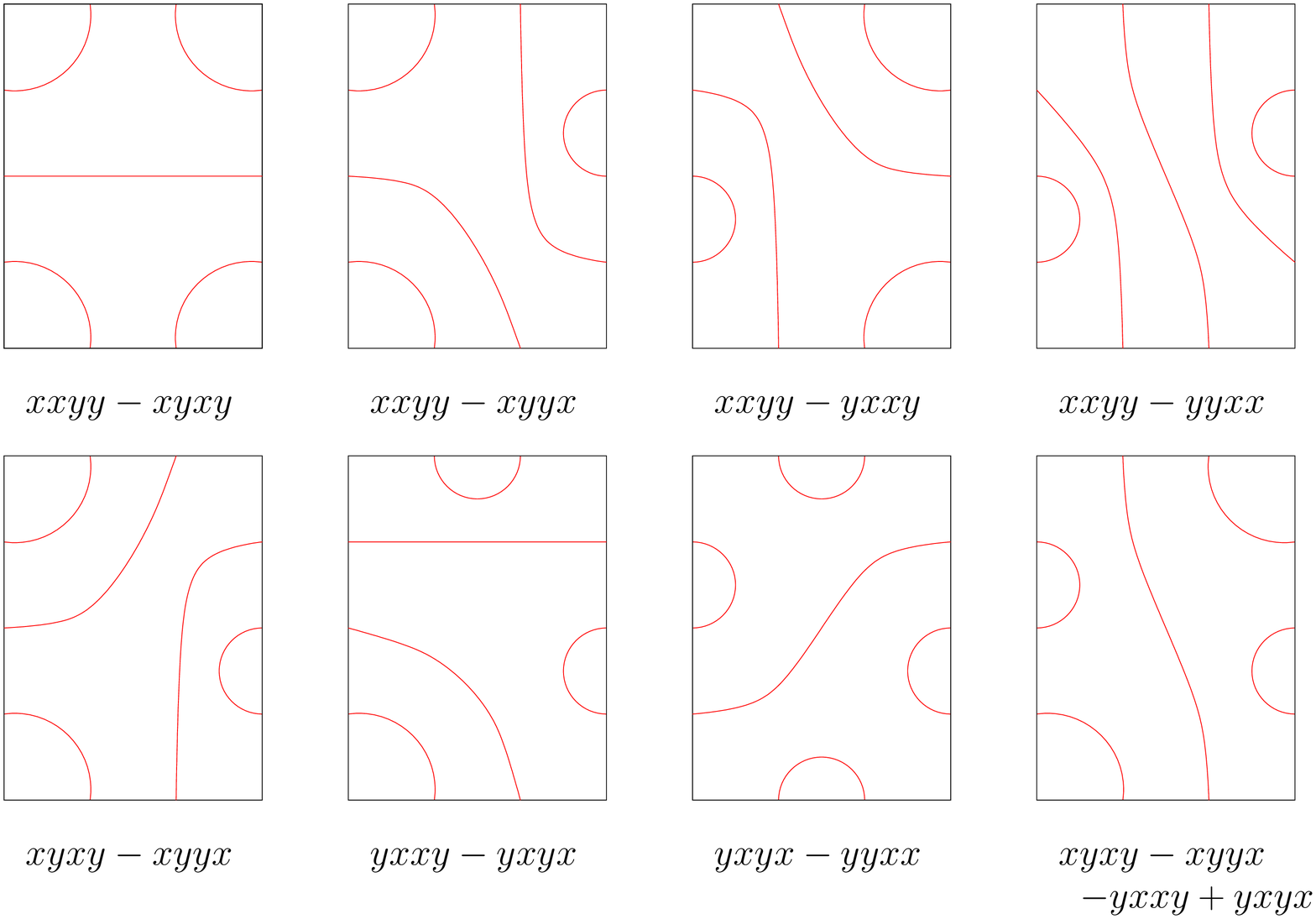}
\caption{Some non-basis chord diagrams with $5$ chords.} \label{fig:25}
\end{figure}

\begin{figure}
\centering
\includegraphics[scale=0.4]{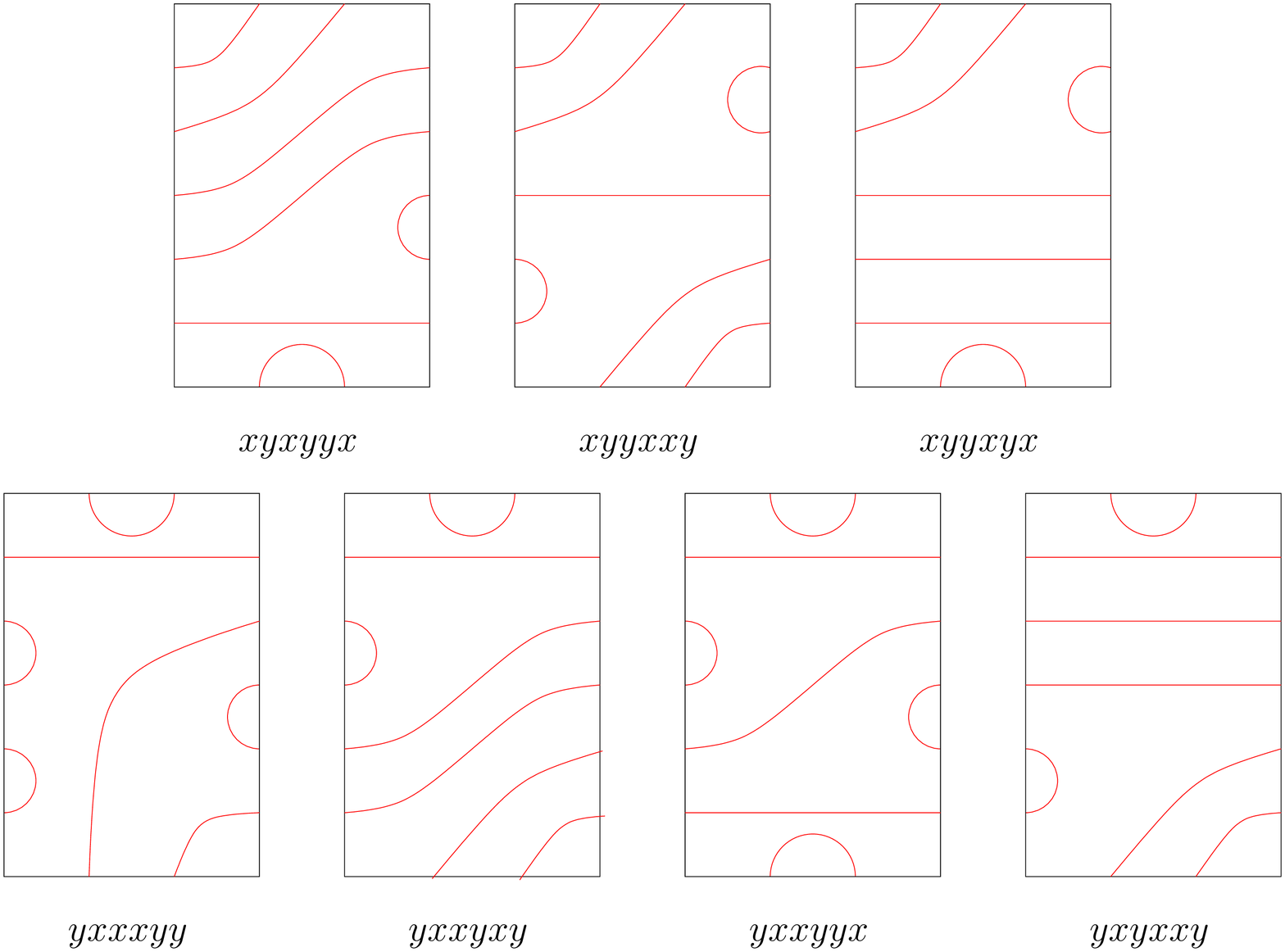}
\caption{Some basis chord diagrams with $7$ chords.} \label{fig:26}
\end{figure}

\begin{figure}
\centering
\includegraphics[scale=0.4]{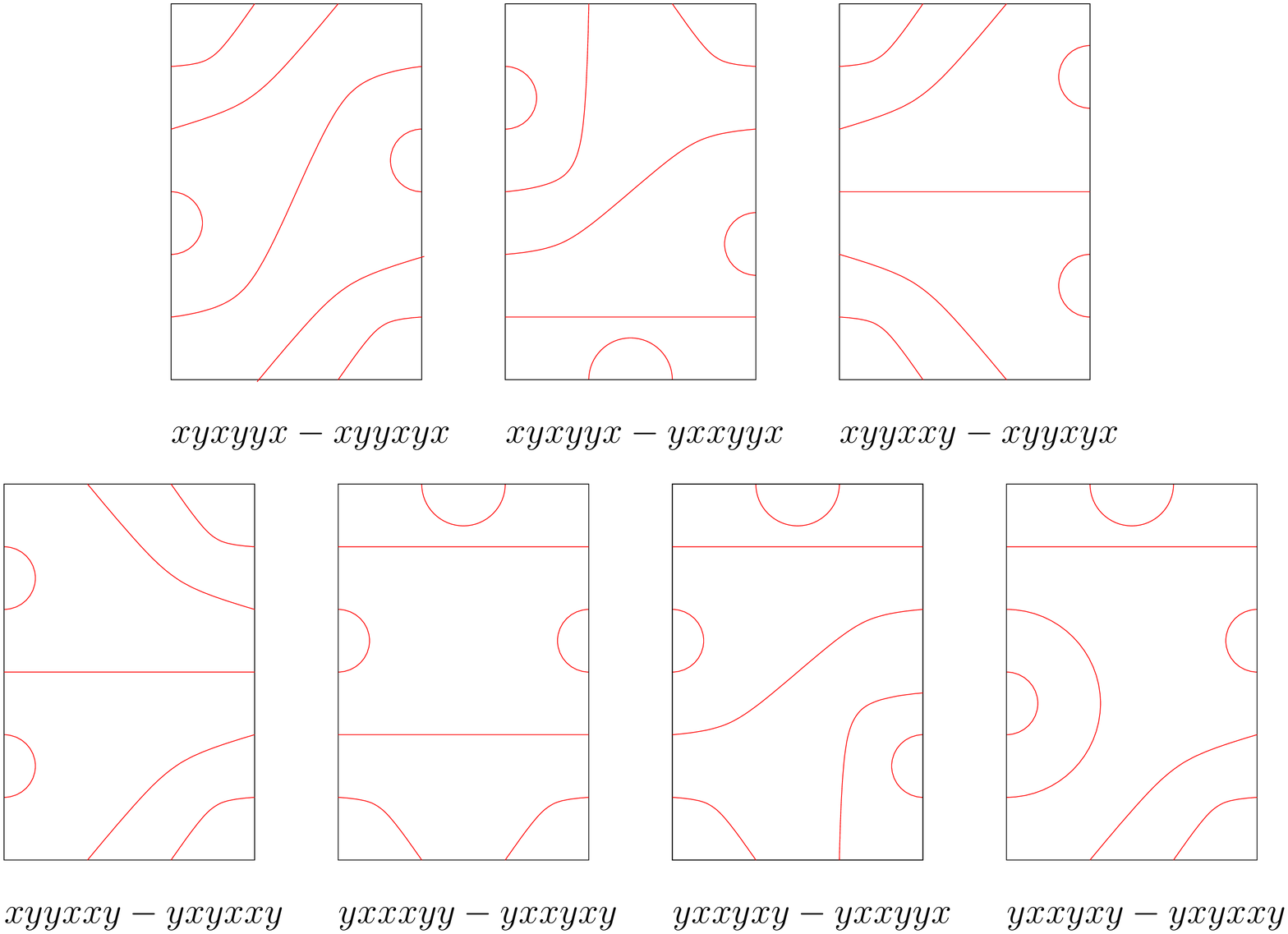}
\caption{Some non-basis chord diagrams with $7$ chords.} \label{fig:27}
\end{figure}

We shall, throughout this paper, consider many chord diagrams. We shall often refer to their suture elements without further explanation. To minimise reader confusion, we provide a menagerie of the chord diagrams we use, and their suture elements, in figures \ref{fig:22}--\ref{fig:27}. Following notation that shall prove useful subsequently, each diagram is drawn in a rectangle, and the basepoint is always taken to be the left marked point on the top side of the rectangle.

\begin{lem}
All of the chord diagrams shown in figures \ref{fig:22}-\ref{fig:27} have a suture element as labelled.
\end{lem}

\begin{Proof}
For the basis elements this follows immediately from the construction algorithm given in \cite{Me09Paper}. For the suture elements which involve two basis elements, we can easily check that the corresponding chord diagram forms a bypass triple with the basis chord diagrams for those two basis elements. Then the fact that coefficients are $\pm 1$ and sum to $0$ gives the suture element shown. 

There is only one suture element shown which involves more than two basis elements, namely $xyxy-xyyx-yxxy+yxyx$. The corresponding diagram forms a bypass triple with the diagrams corresponding to $xyxy-xyyx$ and $yxxy-yxyx$; hence is $xyxy-xyyx \pm(yxxy-yxyx)$ and is the unique (up to sign) suture element starting with $xyxy$ and ending with $yxyx$. As suture elements are closed under multiplication, $(xy-yx)(xy-yx) = xyxy-xyyx-yxxy+yxyx$ is a suture element, hence a suture element of the given diagram.
\end{Proof}

\section{Annuli}
\label{sec:annuli}

\subsection{Notation}

We now consider sutured background $(A, F_A)$ where $A$ is an annulus and $F_A$ consists of two marked points on each boundary component. 

\begin{figure}
\centering
\includegraphics[scale=0.5]{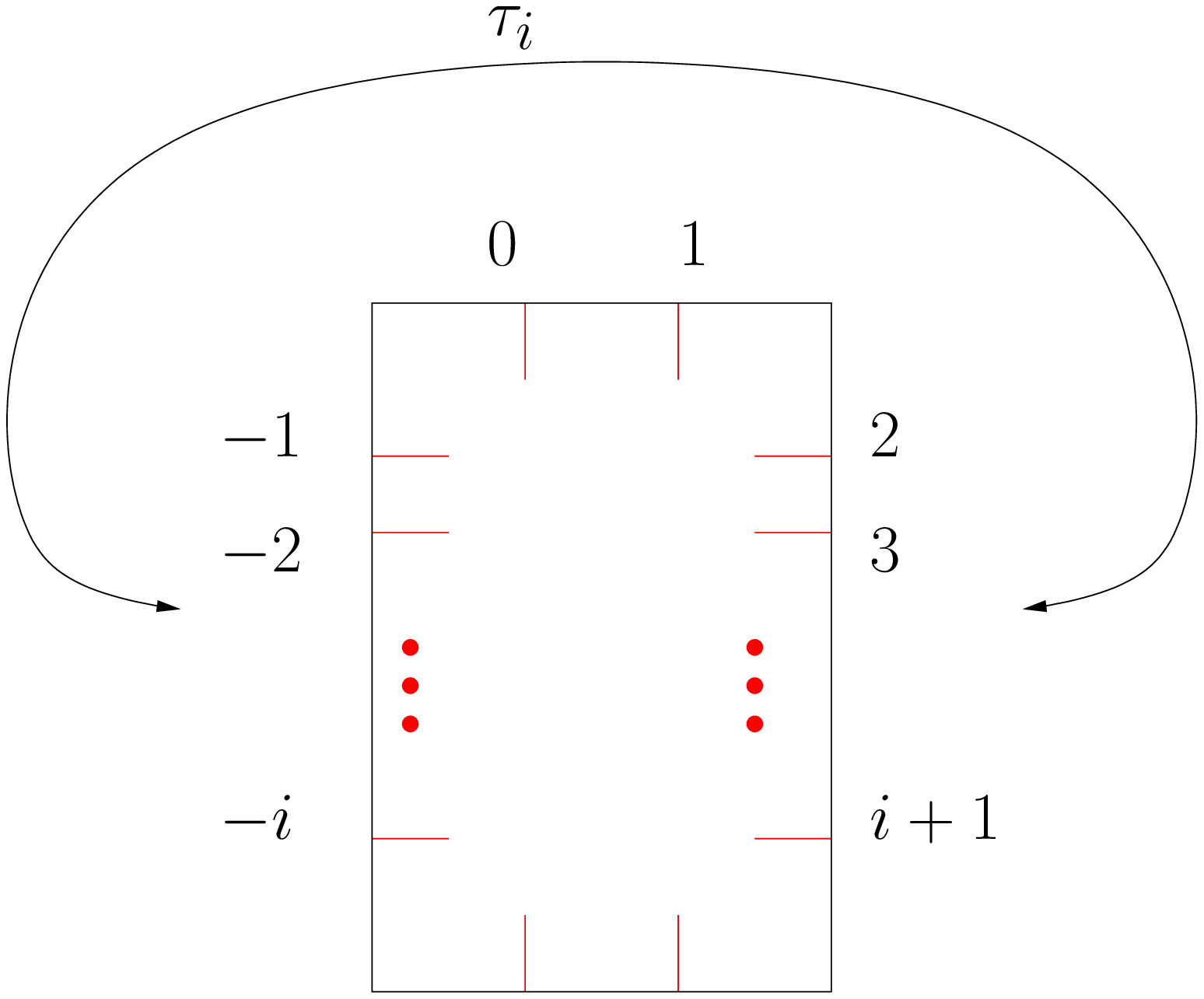}
\caption{Gluing the disc $(D^2, F_{i+2})$ into the annulus $(A,F_A)$ via $\Phi_i$.} \label{fig:14}
\end{figure}

We shall introduce the following notation to keep track of sutures on $(A, F_A)$. Consider $D^2$ as a rectangle; we glue the left and right hand sides together to obtain $(A,F_A)$. The top and bottom of $D^2$ become the two boundary components, hence have two marked points each. We choose the basepoint on $D^2$ to be the leftmost marked point at the top of the rectangle. Let $i$ be the number of marked points on the left and right sides, so there are $2i+4$ marked points in total, and the marked points $(2,3,\ldots,i+1)$ on the right side are glued to $(-1,-2,\ldots,-i)$ on the left side. Let $\tau_i$ denote this gluing, and $\Phi_i$ a map obtained:
\[
\Phi_i: V(D^2, F_{i+2}) \cong V(D^2)_{i+1} \To V(A,F_A).
\]
See figure \ref{fig:14}. This $\Phi_i$ is ambiguous up to sign on each Euler class summand and we will choose signs as we need them.

\subsection{Torsion}
\label{sec:torsion_computation}

In this section we prove the following theorem.
\begin{thm}
\label{thm:torsion_zero}
In any sutured TQFT (satisfying axioms 1--10), if $\Gamma$ is a set of sutures on any $(\Sigma, F)$ with torsion then $c(\Gamma) = \{0\}$.
\end{thm}

\begin{figure}
\centering
\includegraphics[scale=0.3]{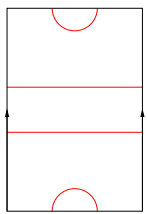}
\caption{Torsion sutures $\Gamma_0$.} \label{fig:13}
\end{figure}

Note that, on any such$(\Sigma, F)$ and $\Gamma$, there is an embedded $(A, F_A, \Gamma_0) \subset (\Sigma, F, \Gamma)$, where $\Gamma_0$ is the set of sutures on $(A, F_A)$ consisting of two boundary-parallel arcs and two parallel closed core curves (figure \ref{fig:13}). We will show that $c(\Gamma_0) = \{0\}$. Then from the inclusion axiom (5') of sutured TQFT, theorem \ref{thm:torsion_zero} follows immediately.

\begin{lem}
\label{lem:torsion_lem1}
We may choose signs on $\Phi_1, \Phi_3$ so that:
\begin{align*}
\Phi_3(xxyy) &= \Phi_1(xy) \\
\Phi_3(xyyx) &= \Phi_1(yx) \\
\Phi_3(yxxy) &= \Phi_1(yx) \\
\Phi_3(yyxx) &= \Phi_1(yx)
\end{align*}
\end{lem}
(Note here we are identifying elements of $V(D^2, F_5)$ with $\F_4$, and elements of $V(D^2, F_3)$ with $\F_2$. The words in $\{x,y\}$ here are basis elements of $\F$, hence correspond to basis chord diagrams.)

\begin{Proof}
Inspect the gluing $\tau_3$ on the basis chord diagrams corresponding to the words $xxyy$, $xyyx$, $yxxy$, $yyxx$ and compare to the gluing $\tau_1$ on the basis chord diagrams corresponding to $xy,yx,yx,yx$. After gluing we obtain isotopic sets of sutures on $(A,F_A)$. Thus $\Phi_3$ applied to these basis elements of $V(D^2, F_5)$ gives the same result as $\Phi_1$ on the corresponding basis elements of $V(D^2, F_3)$, \emph{up to sign}. So we have obtain the four claimed equalities up to sign.

Choose a sign on $\Phi_3$ arbitrarily, then choose a sign on $\Phi_1$ so that $\Phi_3(xxyy) = \Phi_1(xy)$. We claim then that all of the signs are as desired.

We check that $\Phi_3 (xyyx) = \Phi_1(yx)$. Suppose otherwise, so $\Phi_3 (xyyx) = - \Phi_1(yx)$. Then $\Phi_3 (xxyy-xyyx) = \Phi_1(xy+yx)$. But note that $xxyy-xyyx \in \F_4 \cong V(D^2, F_5)$ is a suture element for sutures as shown in figure \ref{fig:15}. Inspecting the gluing of this chord diagram, we see that $\Phi_3(xxyy-xyyx) = \pm \Phi_1(xy-yx)$. Thus $\Phi_1(xy+yx) = \pm \Phi_1 (xy-yx)$, and we have either $2\Phi_1(xy)=0$ or $2\Phi_1(yx)=0$. As $\Phi_1(xy), \Phi_1(yx)$ are nonzero (proposition \ref{prop:nonzero}) and non-torsion (proposition \ref{prop:no_torsion}), this is a contradiction; so $\Phi_3(xyyx)=\Phi_1(yx)$.

Similarly, suppose $\Phi_3(yxxy) = - \Phi_1(yx)$, so $\Phi_3 (xxyy-yxxy) = \Phi_1 (xy+yx)$. Again $xxyy-yxxy$ is a suture element and we find $\Phi_3 (xxyy-yxxy) = \pm \Phi_1(xy-yx)$. Thus $\Phi_1(xy+yx)=\pm \Phi_1(xy-yx)$, impossible as $\Phi_1(xy),\Phi_1(yx)$ are nonzero and non-torsion; so $\Phi_3(yxxy)=\Phi_1(yx)$.

Again, suppose $\Phi_3(yyxx)=-\Phi_1(yx)$, so $\Phi_3(xxyy-yyxx)=\Phi_1(xy+yx)$; inspecting suture elements, we have $\Phi_3(xxyy-yyxx)=\pm \Phi_1(xy-yx)$, so $\Phi_1(xy+yx)=\pm\Phi_1(xy-yx)$ again, a contradiction; hence $\Phi_3(yyxx)=\Phi_1(yx)$.
\end{Proof}

We now define a ``filling in the hole'' gluing map $\Xi$. Consider taking $(A, F_A)$ and gluing a disc to one boundary of it --- the bottom boundary as drawn in our diagrams. Consider the disc glued in to have a set of sutures which is the vacuum. Then we obtain a gluing map
\[
\Xi \; : \; V(A,F_A) \To V(D^2, F_1) = \Z.
\]

\begin{lem}
\label{lem:torsion_lem2}
\[
\Xi \Phi_1 (xy) = \pm 1, \quad \Xi \Phi_1 (yx) = 0, \quad \Xi \Phi_3 (xxyy-xyxy) = 0.
\]
\end{lem}

\begin{Proof}
These are immediate, once we draw the chord diagrams corresponding to the suture elements $xy,yx$ on $(D^2, F_3)$ and $xxyy-xyxy$ on $(D^2, F_5)$, and consider gluing into an annulus and filling in the hole. The chord diagram corresponding to $xxyy-xyxy$ is shown in figure \ref{fig:15}.
\end{Proof}

\begin{lem}
\label{lem:torsion_lem3}
With the choice of $\Phi_1, \Phi_3$ above, $\Phi_3(xyxy) = \Phi_1(xy+yx)$.
\end{lem}

\begin{Proof}
Consider the suture element $xxyy - xyxy$ corresponding to the chord diagram shown in figure \ref{fig:15}: this forms a bypass triple with the suture elements for $xxyy - xyyx$ and $xyyx - xyxy$ respectively. In particular,
\[
xxyy - xyxy = (xxyy - xyyx) + (xyyx - xyxy).
\]
\begin{figure}
\centering
\includegraphics[scale=0.4]{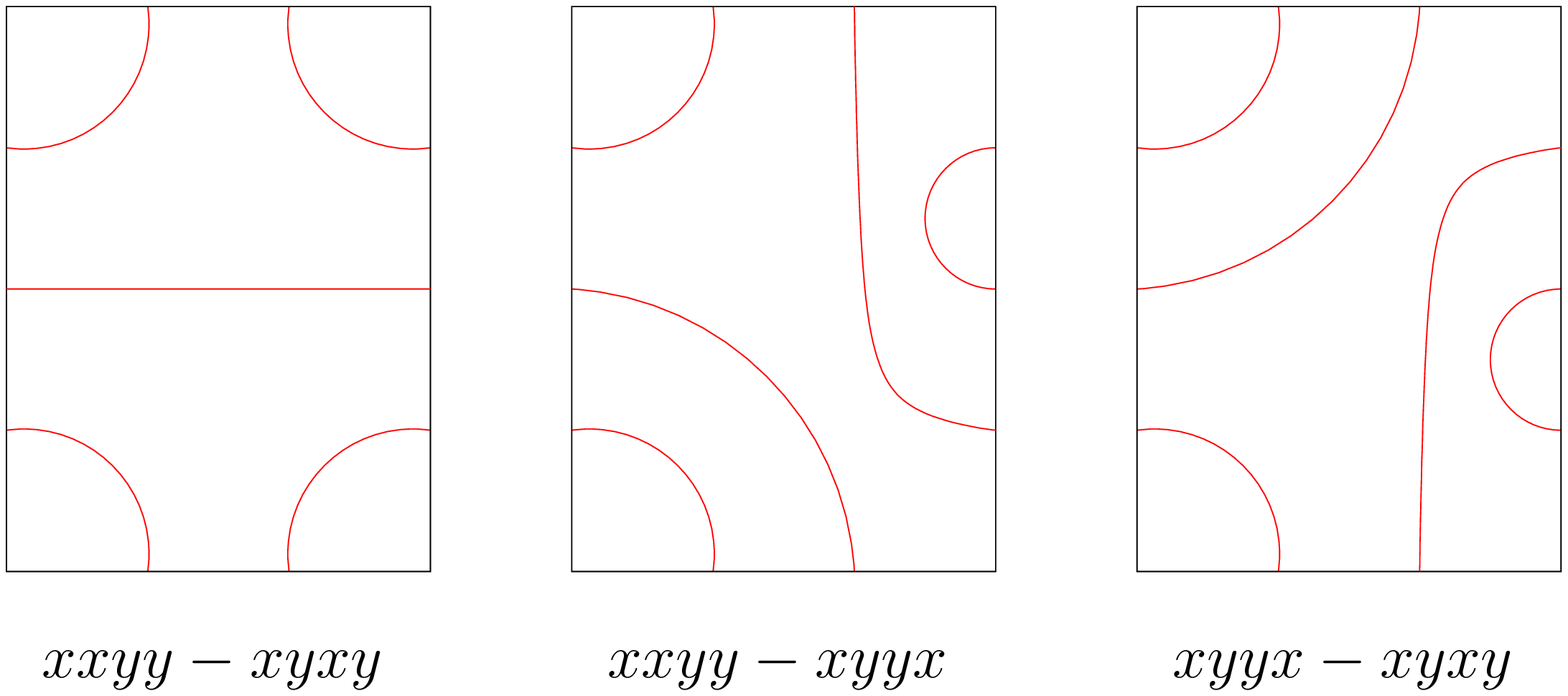}
\caption{Sutures on $(D^2, F_5)$ with suture elements.} \label{fig:15}
\end{figure}
Now apply $\Phi_3$ to these. From lemma \ref{lem:torsion_lem1} we have $\Phi_3(xxyy-xyyx) = \Phi_1(xy-yx)$. Inspecting the gluing of the sutures for $xyyx-xyxy$ we see $\Phi_3(xyyx-xyxy) = \pm \Phi_1(xy)$. So we have
\[
\Phi_3 (xxyy - xyxy) = \Phi_1 (xy-yx) \pm \Phi_1 (xy) = \left\{ \begin{array}{l} \Phi_1 (2xy-yx) \\ \Phi_1 (-yx) \end{array} \right.
\]
Now consider applying $\Xi$. From lemma \ref{lem:torsion_lem2} above we obtain
\[
0 = \Xi \Phi_3 (xxyy - xyxy) = \left\{ \begin{array}{ccc} \Xi \Phi_1 (2xy-yx) & = & \pm 2 \\ \Xi \Phi_1 (-yx) & = & 0 \end{array} \right.
\]
It follows of course that the lower option must be the case, hence $\Phi_3 (xxyy-xyxy) = \Phi_1 (-yx)$. From lemma \ref{lem:torsion_lem1} above we have $\Phi_3(xxyy)=\Phi_1(xy)$, hence $\Phi_3 (xyxy) = \Phi_1 (xy+yx)$ as desired.
\end{Proof}

\begin{lem}
\label{lem:torsion_lem4}
$\Phi_3 (yxyx) = 0$.
\end{lem}

The set of sutures represented by $yxyx$, after gluing into an annulus, has torsion --- in fact, a lot of torsion. So we should obtain zero.

\begin{Proof}
The chord diagram corresponding to $yxyx$ forms a bypass triple with the diagrams given by $yxyx-yyxx$ and $yyxx$:
\[
yxyx = (yxyx-yyxx) + yyxx
\]
Applying $\Phi_3$ to these sets of sutures, from inspection of sutures we have $\Phi_3(yxyx-yyxx) = \pm \Phi_1(yx)$, and from lemma \ref{lem:torsion_lem1} we have $\Phi_3(yyxx)=\Phi_1(yx)$. Thus
\[
\Phi_3(yxyx) = \pm \Phi_1(yx) + \Phi_1(yx) = \left\{ \begin{array}{c} 2 \Phi_1(yx) \\ 0 \end{array} \right.
\]
Now consider the suture element represented by
\[
xyxy - xyyx - yxxy + yxyx,
\]
also shown in figure \ref{fig:25}. Inspecting sutures gives that $\Phi_3 (xyxy-xyyx-yxxy+yxyx)=\pm \Phi_1 (xy-yx)$. Hence
\begin{align*}
\pm \Phi_1 (xy-yx) &= \Phi_3 (xyxy-xyyx-yxxy+yxyx) \\
  &= \Phi_1 (xy+yx) - \Phi_1 (yx) - \Phi_1 (yx) + \left\{ \begin{array}{c} 2 \Phi_1 (yx) \\ 0 \end{array} \right. \\
  &= \left\{ \begin{array}{c} \Phi_1 (xy+yx) \\ \Phi_1(xy - yx) \end{array} \right.
\end{align*}
In the second line here, we used lemma \ref{lem:torsion_lem3} that $\Phi_3(xyxy) = \Phi_1(xy+yx)$, and lemma \ref{lem:torsion_lem1} that $\Phi_3(xyyx) = \Phi_3(yxxy) = \Phi_1(yx)$. As $\Phi_1(xy),\Phi_1(yx)$ are nonzero and non-torsion, the lower option must be the case, and $\Phi_3 (yxyx)=0$.
\end{Proof}

We are now ready to prove theorem \ref{thm:torsion_zero}. As discussed above, a set of sutures has torsion precisely if it contains the set of sutures $\Gamma_0$ of figure \ref{fig:13}; recalling our definition of gluing maps $\Phi_i: V(D^2, F_{i+2}) = \F_{i+1} \To V(A,F_A)$, the sutures $\Gamma_0$ have suture element $\Phi_2 (yxy)$. So we must show $\Phi_2(yxy)=0$.

\begin{Proof}[of theorem \ref{thm:torsion_zero}]
Let $\Gamma_1$ denote the set of sutures on $(A,F_A)$ consisting only of two boundary-parallel arcs enclosing positive discs; so $\Phi_2 (xyy) \in c(\Gamma_1)$. Consider also the chord diagram corresponding to $yxy-xyy$; by inspection, we obtain $\Phi_2 (yxy - xyy) = \pm \Phi_2 (xyy)$, hence $\Phi_2 (yxy) = 0$ or $2 \Phi_2 (xyy)$. Suppose $\Phi_2(yxy) \neq 0$, so $\Phi_2(yxy) = 2 \Phi_2 (xyy)$.

Consider now the gluing map defined by attaching to the annulus $(A, F_A)$ another annulus, to the bottom of our annuli as drawn, with sutures as shown in figure \ref{fig:17}; let this map be $\Upsilon: V(A,F_A) \To V(A,F_A)$.

\begin{figure}
\centering
\includegraphics[scale=0.4]{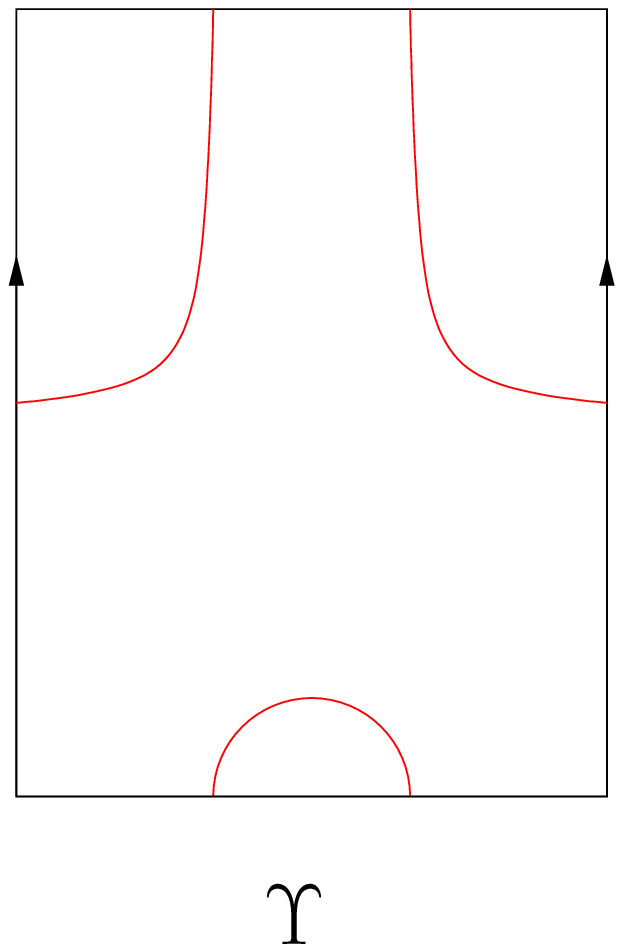}
\caption{The map $\Upsilon$ glues this to the bottom of our annuli as drawn..} \label{fig:17}
\end{figure}

We note by examining sutures that $\Upsilon \Phi_2(xyy) = \pm \Phi_1 (yx)$, and $\Upsilon \Phi_2(yxy) = \pm \Phi_3(yxyx) = 0$ by lemma \ref{lem:torsion_lem4}. But we just deduced that $\Phi_2(yxy) = 2 \Phi_2(xyy)$, hence
\[
0 = \Upsilon \Phi_2(yxy) = 2 \Upsilon \Phi_2(xyy) = \pm 2 \Phi_1(yx).
\]
As $\Phi_1(yx)$ is nonzero and nontorsion, we have a contradiction, and we are done.
\end{Proof}

We remark that the same proof, with essentially the same algebra, can be carried out for different placement of basepoint; in the sequel we will need some details of how the proof works out. We can alternatively choose the basepoint to be the rightmost point on the top of our rectangles for $(D^2, F_n)$. In particular, we then have maps $\Phi_i: V(D^2, F_{i+2}) \To V(A,F_A)$ arising from gluings which identify marked points $(1, \ldots, i)$ to $(-2, \ldots, -i-1)$. The steps involved in the proof above then proceed as follows. For lemma \ref{lem:torsion_lem1} we define $\Phi_3$ arbitrarily and $\Phi_1$ so that $\Phi_3(yyxx) = \Phi_1(yx)$; then we obtain equalities $\Phi_3(xxyy) = \Phi_3(xyyx) = \Phi_3(yxxy) = \Phi_1(xy)$. For lemma \ref{lem:torsion_lem2} we have $\Xi \Phi_1(xy)=0$, $\Xi \Phi_1(yx)= \pm 1$ and $\Xi \Phi_3(yyxx-yxyx)=0$. For lemma \ref{lem:torsion_lem3} we show $\Phi_3(yxyx) = \Phi_1(xy+yx)$: we consider applying $\Phi_3$ to $yyxx-yxyx = (yyxx-yxxy)+(yxxy-yxyx)$ and obtain $\Phi_3 (yyxx-yxyx) = \Phi_1(-xy)$ or $\Phi_1(2yx-xy)$; applying $\Xi$ then shows the first possibility is true, giving $\Phi_3(yxyx)=\Phi_1(xy+yx)$. This is as much as we need.

\subsection{Dehn twists and suture elements}
\label{sec:Dehn_twist_suture_elements}

We now classify suture elements in $V(A,F_A)$. Since we are assuming all axioms \ref{ax:1}--\ref{ax:10}, we have by proposition \ref{prop:no_torsion}, $V(A,F_A) \cong \Z^4$, and $\Phi_1$ as defined above is an isomorphism $V(D^2, F_3) \cong \F_2 \To V(A,F_A)$; and $\Phi_1(xx), \Phi_1(xy), \Phi_1(yx), \Phi_1(yy)$ form a basis for $V(A,F_A)$. We see also that $\Phi_1$ respects euler classes, so that $V(A,F_A)$ splits as a sum $\oplus_e V(A,F_A)^e$ over the possible euler classes $e=-2,0,2$:
\[
\Z^4 \cong V(A,F_A) \cong V(A,F_A)^2 \oplus V(A,F_A)^0 \oplus V(A,F_A)^{-2} \cong \Z \oplus \Z^2 \oplus \Z.
\]

Consider a Dehn twist around the core of $(A,F_A)$. This acts on sets of sutures and can be represented as an inclusion $(A,F_A) \hookrightarrow (A,F_A)$ with prescribed sutures on the intermediate annulus. Hence by axiom (3'), we have a Dehn twist map on sutured TQFT:
\[
\Theta \; : \; V(A,F_A) \To V(A,F_A).
\]
There is sign ambiguity in $\Theta$; we choose a particular representative below. Clearly $\Theta$ respects euler class of suture elements and hence restricts to summands $V(A,F_A)^e \To V(A,F_A)^e$.

\begin{thm}
\label{thm:annulus_suture_elt_classification}
In any sutured TQFT (satisfying axioms \ref{ax:1}--\ref{ax:10}), the nonzero suture elements in $V(A,F_A)$ are precisely as follows.
\begin{enumerate}
\item
The only set of sutures $\Gamma$ on $(A,F_A)$ with $e=-2$ and nonzero suture element consists of two boundary-parallel arcs, each enclosing a negative disc, and $c(\Gamma) = \{\pm 1\}$ in $V(A,F_A)^{-2} \cong \Z$. The Dehn twist map $\Theta$ may be chosen to be the identity on $V(A,F_A)^{-2}$.
\item
The only sets of sutures on $(A,F_A)$ with nonzero suture elements in $V(A,F_A)^0 \cong \Z^2$ consist of the following.
\begin{enumerate}
\item
A closed loop around the core of $(A,F_A)$, and two boundary-parallel arcs; this may be taken to have suture element $\pm(0,1)$.
\item
Two parallel arcs between boundary components, traversing the core of $(A,F_A)$ $n$ times, where $n \in \Z$ (i.e. ``having slope $n/1$''). This may be taken to have suture element $\pm(1,n)$.
\end{enumerate}
With these coordinates, $\Theta$ may be taken to have matrix $\begin{pmatrix} 1 & 0 \\ -1 & 1 \end{pmatrix}$.
\item
The only set of sutures $\Gamma$ with $e=2$ and nonzero suture element consists of two boundary-parallel arcs, each enclosing a positive disc, and $c(\Gamma)= \pm 1$ in $V(A,F_A)^2 \cong \Z$; $\Theta$ may be taken to be the identity on this summand.
\end{enumerate}
\end{thm}

\begin{Proof}
Let $\Gamma$ be a set of sutures on $(A,F_A)$. If $\Gamma$ has a contractible component then clearly $c(\Gamma) = 0$, so assume all closed components are non-contractible, hence homotopic to the core of $A$. If $\Gamma$ has two closed components, then they are parallel and $\Gamma$ has torsion, so from theorem \ref{thm:torsion_zero} $c(\Gamma) = 0$; so assume $\Gamma$ has at most one closed component. Up to Dehn twists then, $\Gamma$ can only take one of four forms shown in figure \ref{fig:18}, corresponding to the cases $e = -2, 0, 0, 2$.

\begin{figure}
\centering
\includegraphics[scale=0.4]{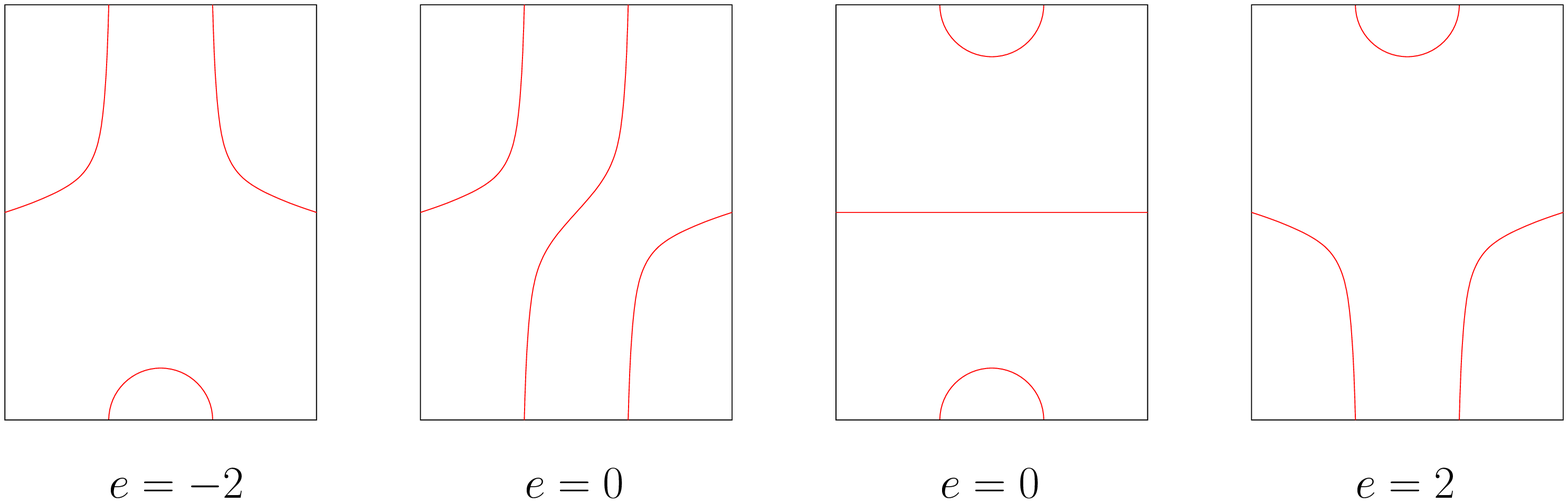}
\caption{Possible sutures on $(A,F_A)$ with nonzero suture elements.} \label{fig:18}
\end{figure}

In the first case $e=-2$ and  $c(\Gamma) = \pm \Phi_1 (xx)$, where $xx \in V(D^2)_2^{-2} \cong \Z$ and $xx$ is a generator for $\Z$, so $c(\Gamma) = \pm 1$. Moreover $\Theta$ clearly takes $c(\Gamma) \to c(\Gamma)$ and hence, adjusting by sign if necessary, we may take $\Theta = 1$. The final case, $e=2$, is similar. We now assume $e=0$.

In this case, the two sutures shown in figure \ref{fig:18} are $\Phi_1 (xy)$ and $\Phi_1 (yx)$; these form a basis for $V(D^2)_2^0$, hence also for $V(A,F_A)^0$, and we denote them $(1,0)$ and $(0,1)$. Now we have by inspecting sutures:
\begin{align*}
\Theta(1,0) &= \Theta \Phi_1 (xy) = \pm \Phi_1 (xy-yx) = \pm(1,-1) \\
\Theta^{-1} (1,0) &= \Theta^{-1} \Phi_1 (xy) = \pm \Phi_3 (xyxy) = \pm \Phi_1(xy+yx) = \pm (1,1), \\
\Theta(0,1) &= \Theta \Phi_1 (yx) = \pm \Phi_1 (yx) = \pm (0,1).
\end{align*}
(In the second line we used lemma \ref{lem:torsion_lem3}.) Altering by a sign if necessary, we set $\Theta (0,1) = (0,1)$. If $\Theta(1,0)=(-1,1)$ then $\Theta(1,1) = \Theta(0,1) + \Theta(1,0) = (0,1) + (-1,1) = (-1,2)$, contradicting $\Theta^{-1}(1,0) = \pm (1,1)$; thus $\Theta(1,0) = (1,-1)$, giving the desired form. As any set of sutures without contractible components or torsion can be obtained by applying Dehn twists to $\pm \Phi_1 (xy) = \pm (0,1)$ or $\Phi_1 (yx) = \pm (1,0)$, the suture elements are precisely $\pm(1,n)$ for $n \in \Z$.
\end{Proof}

By TQFT-inclusion, the above description of $\Theta$ applies any time a Dehn twist is performed on an annulus in which sutures intersect each boundary curve in two points.

\section{Punctured tori}
\label{sec:punctured_tori}

We next consider the sutured background surface $T_1 = (T,F_T)$ where $T$ is a punctured torus and $F_T$ consists of two points on the boundary; the simplest sutured background punctured torus.

\subsection{Gluings of octagons}

\begin{figure}
\centering
\includegraphics[scale=0.4]{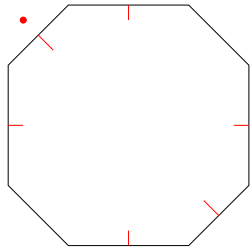}
\caption{The octagon $O_{11}$. Basepoint is marked.} \label{fig:19}
\end{figure}

As with annuli, our strategy is to obtain $(T,F_T)$ by gluing up a sutured background disc, and then use our detailed results for discs. Therefore we consider an sutured background disc in the form of the octagon as shown in figure \ref{fig:19}, which we denote $O_{11}$; it has 6 marked points on the boundary, with one each on left, right, top and bottom sides. The basepoint is placed in the the top left.

We may glue top to bottom and left to right to obtain $(T,F_T)$. The gluing of left to right gives a map $\Psi_{11}: V(O_{11}) \To V(A_1)$, where $A_1 = (A,F_A)$ is the annulus with one pair of marked points on each boundary component; then gluing top to bottom gives a map $\Omega_{1}: V(A_1) \To V(T_1)$. Alternatively, we may first glue top to bottom to obtain a map $\Omega_{11}: V(O_{11}) \To V(A_1)$, and then glue left to right with a map $\Psi_1: V(A_1) \To V(T_1)$. Since an annulus can be obtained in these two distinct ways, write $A_1$ for the annulus obtained from gluing left side to right side, and $B_1$ for the annulus obtained from gluing top to bottom. 

\begin{figure}
\centering
\[
\xymatrix{
O_{33} \ar[dr]^{\Psi_{33}} && O_{31} \ar[dl]_{\Psi_{31}} \ar[dr]^{\Omega_{31}} && O_{11} \ar[dl]_{\Omega_{11}} \ar[dr]^{\Psi_{11}} && O_{13} \ar[dl]_{\Psi_{13}} \ar[dr]^{\Omega_{13}} && O_{33} \ar[dl]_{\Omega_{33}} \\
& A_3 \ar[drrr]_{\Omega_3} && B_1 \ar[dr]^{\Psi_1} && A_1 \ar[dl]_{\Omega_1} && B_3 \ar[dlll]^{\Psi_3} & \\
&&&& T_1 &&&&
}
\]
\caption{Diagram of gluings of octagons into annuli and punctured tori.} \label{fig:20}
\end{figure}
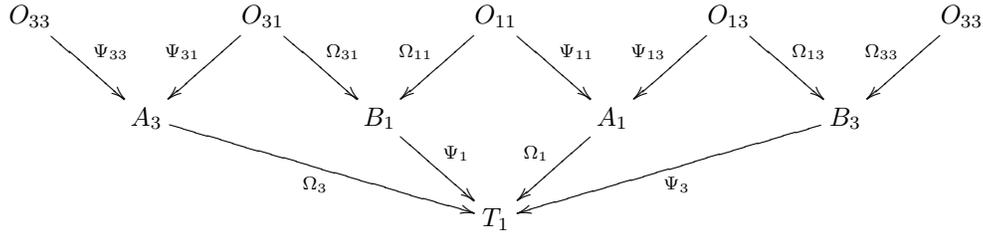

In fact, we will need to consider more general gluings. Let $O_{ij}$ denote the sutured background octagon with $2i+2j+2$ marked points as follows: $i$ points each on top and bottom sides; $j$ points each on left and right sides; a base point on the top-left side; and another point on the bottom-right. Let $\Psi_{ij}$ be the map obtained by gluing the left and right sides of $O_{ij}$, giving an annulus $A_i$ with $i+1$ marked points on each boundary component; then let $\Omega_i$ be the map obtained by gluing top and bottom of the octagon, gluing the annulus $A_i$ into $T_1$. Alternatively, let $\Omega_{ij}$ be the map obtained by gluing the top and bottom of $O_{ij}$, giving an annulus $B_j$ with $j+1$ marked points on each boundary component; then let $\Psi_j$ be the map obtained by gluing left and right of the octagon, gluing the annulus $B_j$ into $T_1$. For these to be sutured background surfaces, we must have $i,j$ both odd. We will consider the four simplest cases $(i,j) = (1,1), (1,3), (3,1), (3,3)$; thus we have the schematic diagram of our gluings in figure \ref{fig:20}.

(Note that $O_{33}$ appears on the left and right of the diagram; we should think of the diagram as ``wrapping around'', so the diagram really contains four squares which ought each to commute, not three.)

Now as $O_{11} \cong (D^2, F_3)$ we have $V(O_{11}) \cong V(D^2, F_3) \cong \Z^4 \cong \Z \oplus \Z^2 \oplus \Z$ over euler classes. By axiom \ref{ax:10}, all of $\Psi_{11}, \Omega_1, \Omega_{11}, \Psi_1$, are isomorphisms. Thus $V(T_1) \cong \Z^4$, and as the isomorphisms preserve euler classes of suture elements, the decomposition $V(T_1) = \oplus_e V(T_1)^e$ is
\[
V(T_1) \cong V(T_1)^{-2} \oplus V(T_1)^0 \oplus V(T_1)^2 \cong \Z \oplus \Z^2 \oplus \Z.
\]

\subsection{Sutures on $T_1$}
\label{sec:sutures_on_T_1}

Consider possible sets of sutures $\Gamma$ on $T_1$. If $\Gamma$ has any contractible component, then $c(\Gamma) = 0$; so assume $\Gamma$ has no contractible component. If $\Gamma$ has two or more closed loop components, then they must be parallel, and $T_1 \backslash \Gamma$ contains an annulus bounded by two sutures, thus has torsion; thus $c(\Gamma) = 0$ again. Thus we may assume $\Gamma$ has at most one closed component. If $\Gamma$ has no closed components, then it must consist of precisely one arc joining the two marked points on the boundary. To produce a valid set of sutures, $\Gamma$ must separate $T_1$; thus the arc is boundary parallel and the Euler class is $e= \pm 2$. If $\Gamma$ contains one closed component, then $\Gamma$ must have two components: the closed loop and an arc connecting the two marked points. If the closed loop is boundary parallel, then so is the arc, and $e=0$. If the closed loop is not boundary parallel then it has some slope; in order to produce a valid set of sutures, the arc then cannot be boundary parallel, but must have the same slope as the closed loop, and $e=0$. In this case note that Dehn twists around the boundary will produce infinitely many distinct sets of sutures all with the same slope.

This gives a complete classification of sets of sutures $\Gamma$ for which $c(\Gamma)$ might be nonzero:
\begin{enumerate}
\item
$e=-2$. One arc, boundary parallel, enclosing a negative disc.
\item
$e=0$.
\begin{enumerate}
\item
One arc and one closed loop, both boundary parallel. The arc may enclose a positive or a negative region.
\item
One arc and one closed loop, both of some slope $q/p$. There is a countably infinite set of sutures with slope $q/p$, arising from Dehn twists around the boundary.
\end{enumerate}
\item
$e=2$. One arc, boundary parallel, enclosing a positive disc.
\end{enumerate}

Suture elements in the cases $e=\pm 2$ are simple enough. For $e=-2$, the suture $\Gamma$ with one boundary parallel arc enclosing a negative disc is obtained from the chord diagram with suture element $xx$ on $O_{11} \cong (D^2, F_3)$, after gluing via $\Psi_1 \circ \Omega_{11}$ or $\Omega_1 \circ \Psi_{11}$. Thus $\Psi_1 \Omega_{11} (xx)$ or $\Omega_1 \Psi_{11} (xx)$ generates $V(T_1)^{-2} \cong \Z$. Similarly for $e=2$, $V(T_1)^2 \cong \Z$ is generated by $\Psi_1 \Omega_{11} (yy)$ or $\Omega_1 \Psi_{11} (yy)$.

Thus, we need henceforth only consider sutures and summands with euler class $0$.

\subsection{Coherent signs for gluings}
\label{sec:coherent_signs_for_gluings}

Note that all the maps $\Omega_{ij}, \Psi_{ij}$, etc., can be adjusted up to sign. We will now choose signs for all of them so that figure \ref{fig:20} commutes. First note that the maps have certain clear relations: various sets of sutures on $O_{ij}$, after gluing, correspond to the same sets of sutures on annuli and tori. For instance, $\Omega_{33}(yxyxxy)$ and $\Omega_{13}(xyxy)$ refer to the same set of sutures on $B_3$, so are equal up to sign. We define the maps as follows. Begin with some arbitrary choices, and some natural ones.
\begin{itemize}
\item 
Choose signs for $\Omega_{33}, \Psi_{33}, \Omega_{11}, \Psi_{11}$ arbitrarily.
\item
Choose a sign for $\Psi_{31}$ so that $\Psi_{31}(yxyx) = \Psi_{33}(yxxyxy)$.
\item
Choose a sign for $\Omega_{13}$ so that $\Omega_{13}(xyxy) = \Omega_{33}(yxyxxy)$.
\end{itemize}

We now note that gluing with $\Psi_{11}$ or $\Psi_{13}$, we have glued left to right of our octagon, giving the annulus $A_1$ with two marked points on each boundary component. Moreover, the gluing of the disc, and position of basepoint, is identical to that considered in section \ref{sec:torsion_computation} above to give annuli; $\Psi_{11}$ is like $\Phi_1$ and $\Psi_{13}$ is like $\Phi_3$. So we may choose $\Psi_{13}$ so that they obey the same relations: in particular,
\begin{itemize}
\item
Choose a sign for $\Psi_{13}$ so that $\Psi_{13}(xyxy) = \Psi_{11}(xy+yx)$ and $\Psi_{13}(xxyy) = \Psi_{11}(xy)$.
\end{itemize}

A similar situation occurs with $\Omega_{11}$ and $\Omega_{31}$; again we are gluing two opposite sides of the octagon to give the annulus $B_1$; however now the basepoint is in a different position to that in section \ref{sec:torsion_computation}. But we can apply the remark at the end of section \ref{sec:torsion_computation}, noting that $\Omega_{11}$ is like $\Phi_1$ there and $\Omega_{31}$ is like $\Phi_3$, and we choose $\Omega_{31}$ so that they obey the same relations; in particular including the following relations.
\begin{itemize}
\item
Choose a sign for $\Omega_{31}$ so that $\Omega_{31}(yxyx)=\Omega_{11}(xy+yx)$ and $\Omega_{31}(yxxy) = \Omega_{11}(xy)$.
\end{itemize}

It remains to define the maps $\Omega_i$ and $\Psi_i$ from annuli to tori. Gluing an octagon's sides in either order ought to give the same result; there ought to be four commutation relations arising from each of $O_{11}, O_{13}, O_{31}, O_{33}$; corresponding to the four squares in figure \ref{fig:20}. We define one map arbitrarily; the other three are then defined so that three of these commutation relations are satisfied; and then we show that the fourth also holds.
\begin{itemize}
\item Choose a sign for $\Omega_1$ arbitrarily.
\end{itemize}
Considering the sets of sutures corresponding to $xy, yx$ and $xy-yx$ on $O_{11}$, we have immediately that
\[
\Omega_1 \Psi_{11} (xy) = \pm \Psi_1 \Omega_{11} (xy), \quad \Omega_1 \Psi_{11} (yx) = \pm \Psi_1 \Omega_{11} (yx), \quad \Omega_1 \Psi_{11} (xy-yx) = \pm \Psi_1 \Omega_{11} (xy-yx).
\]
From linearity of the maps, and the fact that all the suture elements involved are nonzero (proposition \ref{prop:nonzero}) and nontorsion, it follows immediately that all three signs must be the same; as these suture elements span $V(D^2, F_3)^0$, $\Psi_1 \Omega_{11} = \pm \Omega_1 \Psi_{11}$ as maps (on euler class $0$ summands).
\begin{itemize}
\item Choose a sign for $\Psi_1$ so that the gluing square for $O_{11}$ commutes, i.e. $\Psi_1 \Omega_{11} = \Omega_1 \Psi_{11}$.
\end{itemize}
A similar argument applies to the other octagons $O_{ij}$; on each basis element (word) $w$ in $V(D^2, F_{i+j+1})^0$, $\Psi_j \Omega_{ij} (w) = \pm \Omega_i \Psi_{ij} (w)$; for two words $w,w'$ related by an elementary move $\Psi_j \Omega_{ij}(w-w') = \pm \Omega_i \Psi_{ij} (w-w')$; hence the signs are the same for words $w, w'$ related by an elementary move; and thus the signs are the same for all words, and $\Psi_i \Omega_{ij} = \pm \Omega_i \Psi_{ij}$ on the euler class $0$ summand $V(D^2, F_{i+j+1})^0$.
\begin{itemize}
\item Choose a sign for $\Omega_3$ so that the gluing square for $O_{31}$ commutes, $\Omega_3 \Psi_{31} = \Psi_1 \Omega_{31}$.
\item Choose a sign for $\Psi_3$ so that the gluing square for $O_{13}$ commutes, $\Omega_1 \Psi_{13} = \Psi_3 \Omega_{13}$.
\end{itemize}

Now three of the gluing squares commute; we check that the gluing square for $O_{33}$ commutes.
\begin{lem}
with these choices of sign, $\Omega_3 \Psi_{33} = \Psi_3 \Omega_{33}$. Hence figure \ref{fig:20} completely commutes.
\end{lem}

\begin{Proof}
From the preceding argument these two maps agree up to sign; it is therefore sufficient to check that they agree on one basis element. We will show that they agree on $yxxxyy$. We chase the diagram around, using our choices of signs. We will show that paths from the following elements agree.
\[
\xymatrix{
yxxxyy \ar[dr]^{\Psi_{33}} && yxxy \ar[dl]_{\Psi_{31}} \ar[dr]^{\Omega_{31}} && xy \ar[dl]_{\Omega_{11}} \ar[dr]^{\Psi_{11}} && xxyy \ar[dl]_{\Psi_{13}} \ar[dr]^{\Omega_{13}} && yxxxyy \ar[dl]_{\Omega_{33}} \\
& \cdot \ar[drrr]_{\Omega_3} && \cdot \ar[dr]^{\Psi_1} && \cdot \ar[dl]_{\Omega_1} && \cdot \ar[dlll]^{\Psi_3} & \\
&&&& \cdot &&&&
}
\]

First, $\Psi_{33} (yxxxyy) = \Psi_{31}(yxxy)$. This is true by inspection of sutures, up to sign; suppose not, so $\Psi_{33} (yxxxyy) = - \Psi_{31} (yxxy)$. We chose signs so that $\Psi_{33}(yxxyxy) = \Psi_{31}(yxyx)$; so $\Psi_{33}(yxxxyy-yxxyxy) = \Psi_{31}(-yxxy-yxyx)$. On the other hand, by inspection of sutures, $\Psi_{33}(yxxxyy-yxxyxy)= \pm \Psi_{31}(yxxy-yxyx)$. Thus $\Psi_{31}(-yxxy-yxyx) = \pm \Psi_{31}(yxxy-yxyx)$. As these sets of sutures are all non-isolating, by proposition \ref{prop:nonzero} they have nonzero (and nontorsion) suture elements, so we have a contradiction. Thus $\Psi_{33}(yxxxyy) = \Psi_{31}(yxxy)$. Hence $\Omega_3 \Psi_{33} (yxxxyy) = \Omega_3 \Psi_{31}(yxxy)$; since the gluing square for $O_{31}$ commutes, this also equals $\Psi_1 \Omega_{31} (yxxy)$.

Second, $\Omega_{31}(yxxy) = \Omega_{11}(xy)$. This is true from our definition of $\Omega_{31}$, and the remark at the end of section \ref{sec:torsion_computation}. Thus $\Psi_1 \Omega_{31}(yxxy) = \Psi_1 \Omega_{11}(xy)$; by commutativity of the gluing square for $O_{11}$ this is also $\Omega_1 \Psi_{11}(xy)$.

Third, $\Psi_{11}(xy) = \Psi_{13}(xxyy)$, by from definition of $\Psi_{13}$ (which used lemma \ref{lem:torsion_lem1}). Thus $\Omega_1 \Psi_{11}(xy) = \Omega_1 \Psi_{13} (xxyy) = \Psi_3 \Omega_{13} (xxyy)$, by the gluing square for $O_{13}$.

Finally, $\Omega_{13} (xxyy) = \Omega_{33} (yxxxyy)$. This is similar to our first claim. It follows from our choice of $\Omega_{13}$ so that $\Omega_{13}(xyxy) = \Omega_{33}(yxyxxy)$, together with the equalities up to sign from suture inspection, $\Omega_{13}(xxyy) = \pm \Omega_{33}(yxxxyy)$ and $\Omega_{13}(xxyy-xyxy) = \pm \Omega_{33} (yxxxyy-yxyxxy)$ (and the fact that these sutures are all non-isolating, hence have nonzero suture elements). Then $\Psi_3 \Omega_{13}(xxyy) = \Psi_3 \Omega_{33} (yxxxyy)$. 

We have now come full circle: we have $\Omega_3 \Psi_{33} (yxxxyy) = \Psi_3 \Omega_{33} (yxxxyy)$, so that $\Omega_3 \Psi_{33} = \Psi_3 \Omega_{33}$ agree on one, hence all of a generating set.
\end{Proof}

\subsection{Boundary parallel sutures}

It follows from our definitions and sign choices that the paths from the following elements commute also:
\[
\xymatrix{
yxxyxy \ar[dr]^{\Psi_{33}} && yxyx \ar[dl]_{\Psi_{31}} \ar[dr]^{\Omega_{31}} && xy+yx \ar[dl]_{\Omega_{11}} \ar[dr]^{\Psi_{11}} && xyxy \ar[dl]_{\Psi_{13}} \ar[dr]^{\Omega_{13}} && yxyxxy \ar[dl]_{\Omega_{33}} \\
& \cdot \ar[drrr]_{\Omega_3} && \cdot \ar[dr]^{\Psi_1} && \cdot \ar[dl]_{\Omega_1} && \cdot \ar[dlll]^{\Psi_3} & \\
&&&& \cdot &&&&
}
\]
Thus $\Omega_{3} \Psi_{33} = \Psi_3 \Omega_{33}$ agree on $yxyxxy$ and $yxxyxy$. It follows that the set of sutures on $T_1$ represented by $\Omega_3 \Psi_{33} (yxxyxy-yxyxxy)$ has suture element $0$.

\begin{thm}
\label{thm:boundary_parallel_sutures}
Let $\Gamma$ be a set of sutures on $T_1$ consisting of a boundary-parallel arc and a boundary-parallel closed loop. Then $c(\Gamma)=0$.
\end{thm}

\begin{Proof}
The boundary-parallel arc either encloses a positive or negative region. If it encloses a positive region then its suture element is given by $\Omega_3 \Psi_{33} (yxxyxy-yxyxxy)$, which is $0$ from above. If it encloses a negative region then its suture element is given by $\Omega_3 \Psi_{33}(xyxyyx-xyyxyx)$; we now show this is also $0$. 

First consider $\Omega_{33} (xyyxyx)$. We note that $\Omega_{33}(xyyxyx-xyyxxy)=0$, since the corresponding set of sutures on the annulus has a contractible loop. Also $\Omega_{33}(xyyxxy-yxyxxy) = 0$, for the same reason. Thus $\Omega_{33}(xyyxyx)=\Omega_{33}(xyyxxy)=\Omega_{33}(yxyxxy)$.

Similarly, consider $\Psi_{33}(xyxyyx)$. We now note that $\Psi_{33}(xyxyyx-yxxyyx)=0$ and $\Psi_{33}(yxxyyx-yxxyxy)=0$, since both are sutures on the annulus with contractible loops. Thus $\Psi_{33}(xyxyyx)=\Psi_{33}(yxxyyx)=\Psi_{33}(yxxyxy)$.

We just showed above that $\Omega_3 \Psi_{33}(yxxyxy-yxyxxy)=0$, hence
\[
0 = \Omega_3 \Psi_{33} (yxxyxy) - \Psi_3 \Omega_{33} (yxyxxy) = \Omega_3 \Psi_{33} (xyxyyx) - \Psi_3 \Omega_{33} (xyyxyx) = \Omega_3 \Psi_{33} (xyxyyx - xyyxyx)
\]
as desired.
\end{Proof}

\begin{cor}
\label{cor:isolated_punctured_torus}
Suppose $\Gamma$ is a set of sutures on a sutured background surface $(\Sigma, F)$, and $\Sigma \backslash \Gamma$ has a component which is a punctured torus. Then $c(\Gamma)=0$.
\end{cor}

\begin{Proof}
There exists an embedded punctured torus $T_1$ with sutures $\Gamma'$ of the type considered in theorem \ref{thm:boundary_parallel_sutures}, embedded in $\Sigma$ with sutures $\Gamma$; hence $c(\Gamma')=0 \in V(T_1)$. By the inclusion properties of axioms ($3',5'$), $c(\Gamma)=0$ also.
\end{Proof}

\subsection{Isolation}
\label{sec:isolation}

We can now prove theorem \ref{thm:STQFT_isolating_vanishing}: in a sutured TQFT satisfying axioms \ref{ax:1}--\ref{ax:10}, an isolating set of sutures on any sutured background surface $(\Sigma, F)$ has suture element $0$.

\begin{Proof}[of theorem \ref{thm:STQFT_isolating_vanishing}]
We apply the argument of Massot's theorem 16 in \cite{Massot09}. Suppose there is an isolated component of $\Sigma \backslash \Gamma$ with genus $g$ and $n$ boundary components. Massot shows that $\Gamma$ forms a bypass triple with two other sets of sutures $\Gamma'$ and $\Gamma''$, such that each of $\Gamma$ and $\Gamma'$ have isolated components which are topologically simpler --- an annulus, or punctured torus, or a surface of genus $<g$, or with $<n$ boundary components). Hence, by the bypass relation, if $c(\Gamma) = 0$ whenever $\Gamma$ has an isolated annulus or punctured torus, then $c(\Gamma)=0$ whenever $\Gamma$ is isolating. 

We proved in theorem \ref{thm:torsion_zero} that for any sutures $\Gamma$ with an isolated annulus (i.e. torsion), $c(\Gamma)= 0$. And we proved in corollary \ref{cor:isolated_punctured_torus} that for sutures $\Gamma$ with an isolated punctured torus, $c(\Gamma)=0$. This immediately gives the theorem.
\end{Proof}

\subsection{Classification}
\label{sec:classification}

We now classify suture elements in $V(T_1)$. Recall from section \ref{sec:sutures_on_T_1} that there remains only Euler class $0$ to consider: $V(T_1)^{\pm 2} \cong \Z$, generated by $\Psi_1 \Omega_{11}$ or $\Omega_1 \Psi_{11}$ applied to $xx$ or $yy$ in $V(D^2, F_3)^{\pm 2}$. 

In euler class $0$, we have $V(T_1)^0 \cong \Z^2$, with a basis given by applying $\Psi_1 \Omega_{11} = \Omega_1 \Psi_{11}$ to $xy$ and $yx$ (these maps are equal by our choices of signs in section \ref{sec:coherent_signs_for_gluings}, and isomorphisms by axiom \ref{ax:10}). We can write $\Psi_1 \Omega_{11} (yx) = (1,0)$ and $\Psi_1 \Omega_{11} (xy) = (0,1)$ in coordinates on $V(T_1)^0$.

Recall from section \ref{sec:sutures_on_T_1} that in euler class $0$, there are two classes of sutures to consider, (a) and (b). In the boundary-parallel case (a) we now know from theorem \ref{thm:boundary_parallel_sutures} that the suture element is $0$. In case (b), we have a closed loop and an arc, both of slope $q/p$; Dehn twisting around the boundary gives distinct sets of sutures. A Dehn twist around the boundary can clearly be achieved by including $T_1$ into a larger $T_1$, with prescribed sutures on an intermediate annulus, being a neighbourhood of the boundary. We show the effect of this Dehn twist, i.e. this inclusion of punctured tori, on sutured TQFT is trivial.
\begin{lem}
\label{lem:boundary_Dehn_twist_1}
We may choose a sign for the map $\Theta: V(T_1)^0 \To V(T_1)^0$, induced by the inclusion $T_1 \hookrightarrow T_1$ corresponding to a Dehn twist about $\partial T_1$, so that $\Theta$ is the identity.
\end{lem}

\begin{Proof}
If $\Gamma$ is a set of sutures on $T_1$ consisting of a loop and arc of slope $q/p$, and $\Gamma'$ is the set of sutures obtained from $\Gamma$ by a Dehn twist about the boundary, then consider bypass surgery along an attaching arc $\delta$ near $\partial T_1$ on $\Gamma'$ as shown in figure \ref{fig:21}: performing bypass surgery in one direction gives an isolating set of sutures, and in the other direction gives $\Gamma$. Hence $c(\Gamma') = \pm 0 \pm c(\Gamma) = \pm c(\Gamma)$.

\begin{figure}
\centering
\includegraphics[scale=0.4]{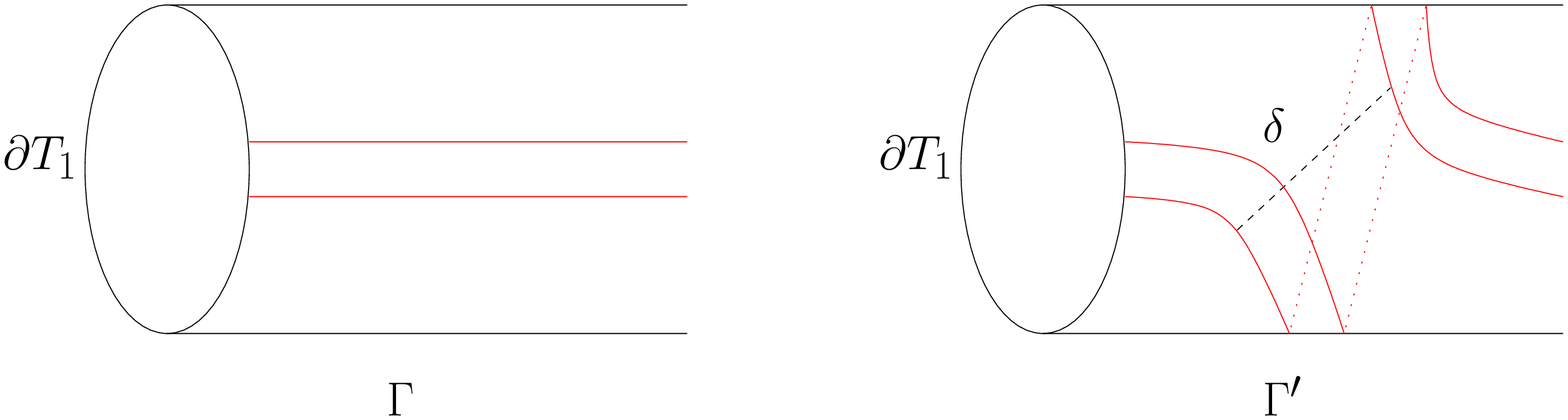}
\caption{Sutures $\Gamma$, and $\Gamma'$ after a Dehn twist around $\partial T_1$.} \label{fig:21}
\end{figure}

Now we have
\begin{align*}
\Theta \Psi_1 \Omega_{11} (xy) &= \pm \Psi_1 \Omega_{11} (xy),\\
\Theta \Psi_1 \Omega_{11} (yx) &= \pm \Psi_1 \Omega_{11} (yx),\\
\Theta \Psi_1 \Omega_{11} (xy-yx) &= \pm \Psi_1 \Omega_{11} (xy-yx).
\end{align*}
Since these three sets of sutures on $T_1$ are non-isolating, they have nonzero (and nontorsion) suture elements. It follows that all the signs are the same; we choose a sign for $\Theta$ so they are all $+$. As $\Psi_1 \Omega_{11} (xy)$, $\Psi_1 \Omega_{11} (yx)$ form a basis for $V(T_1)^0$, then $\Theta = 1$.
\end{Proof}

Alternatively, lemma \ref{lem:boundary_Dehn_twist_1} follows from the computation of the Dehn twist on the annulus in theorem \ref{thm:annulus_suture_elt_classification}, and including it into the punctured torus. Including the sutures on the annulus consisting of a core curve and two boundary-parallel arcs into $T_1$ gives isolating sutures, so maps to zero in $V(T_1)^0$, by theorem \ref{thm:STQFT_isolating_vanishing}.

It now follows that the nonzero suture elements in $V(T_1)^0$ only arise from sutures $\Gamma_{q/p}$ consisting of a loop and arc of a particular slope $q/p$; and $c(\Gamma_{q/p})$ depends only on the slope $q/p$, not on Dehn twists about the boundary. To make a precise statement, note that the slope of a curve is measured relative to some choice of basis for the first homology $H_1(T_1)$. Choose a basis for $H_1(T_1) \cong \Z^2$ so that the loops representing $(0,1)$ and $(1,0)$ in $H_1$ identify with the standard Cartesian coordinate directions in our drawings of octagons, glued into punctured tori; then denote a curve $\pm(p,q) \in H_1(T_1)$ to have slope $q/p$, so this agrees with our standard notion of slope in the Cartesian plane.

Thus, after gluing up the chord diagram with suture element $xy$ to obtain a set of sutures on $T_1$, we have the suture element $\Psi_1 \Omega_{11} (xy)$, which we have defined to have coordinates $(0,1) \in V(T_1)^0$; and the sutures represent $\pm(0,1) \in H_1(T_1)$ and have slope $\infty = 1/0$; hence $c(\Gamma_{1/0}) = \pm(0,1)$. Similarly, the sutures defined by gluing up the chord diagram for $yx$ have suture element $\Psi_1 \Omega_{11} (yx) = (1,0) \in V(T_1)^0$; and the sutures themselves represent $\pm(1,0) \in H_1(T_1)$, with slope $0/1$; so $c(\Gamma_{0/1}) = \pm(1,0)$. Likewise, the sutures given by gluing up the chord diagram $xy-yx$ have suture element $\Psi_1 \Omega_{11}(xy-yx) = (-1,1) \in V(T_1)^0$, and the sutures themselves represent $\pm(-1,1) \in H_1(T_1)$, with slope $-1/1$, so $c(\Gamma_{-1/1} = \pm(-1,1)$.

We then see a correspondence, at least for these three sets of sutures, between the suture elements in $V(T_1)^0$ of sutures, and their homology classes in $H_1(T_1)$. We will show this holds for all sets of sutures $\Gamma_{q/p}$ described above; suture elements, and the homology class of sutures, have the same coordinates. This gives a complete classification result
\begin{thm}
\label{thm:punctured_torus_suture_elt_classification}
The nonzero suture elements in $V(T_1)$ are precisely as follows.
\begin{enumerate}
\item
The only set of sutures $\Gamma$ with $e=-2$ and $c(\Gamma) \neq 0$ consists of one boundary-parallel arc enclosing a negative disc, and $c(\Gamma) = \{\pm 1\}$ in $V(T_1)^{-2} \cong \Z$.
\item
The only sets of sutures with $e=0$ and nonzero suture element are $\Gamma_{q/p}$, consisting of an arc and a closed loop of the same slope $q/p$. Choosing bases for $H_1(T_1)$ and $V(T_1)^0$ as above, $c(\Gamma_{q/p}) = \pm(p,q)$.
\item
The only set of sutures $\Gamma$ with $e=2$ and $c(\Gamma) \neq 0$ consists of one boundary-parallel arc enclosing a positive disc, and $c(\Gamma) = \{\pm 1\}$ in $V(T_1)^2 \cong \Z$.
\end{enumerate}
\end{thm}

In \cite{Conway}, Conway introduced some definitions which are useful in the present context, as also mentioned in \cite{Me10_Sutured_TQFT}. An element $v$ of an abelian group is a \emph{strict element}, and the pair $\{\pm v\}$ is a \emph{lax element} --- when we lazily decline to care about signs. A basis $(v,w)$ of $\Z^2$ is a \emph{strict basis}, and the pair $(\pm v, \pm w)$ of lax elements is a \emph{lax basis}. A \emph{strict superbasis} is a triple $(u,v,w)$ of strict elements such that $(u,v)$ is a (strict) basis and $u+v+w=0$ (hence also $(v,w)$ and $(w,u)$ are bases). A \emph{lax superbasis} is a triple $(\pm u, \pm v, \pm w)$ such that $(u,v,w)$ is a strict superbasis. A lax superbasis $(\pm u, \pm v, \pm w)$ naturally restricts to three distinct lax bases: $(\pm u, \pm v)$, $(\pm v, \pm w)$ and $(\pm w, \pm u)$. A lax basis $(\pm v, \pm w)$ naturally includes into two distinct lax superbases: $(\pm (v+w), \pm v, \pm w)$ and $(\pm (v-w), \pm v, \pm w)$.

This terminology is useful here. A slope $q/p = (-q)/(-p)$ can be considered a lax element $\pm(p,q) \in \Z^2$. A bypass surgery on $\Gamma_{q/p}$, consisting of an arc and loop of slope $q/p$, is either trivial or produces $\Gamma_{s/r}$, where $ps-qr = \pm 1$, i.e. $\{(p,q),(r,s)\}$ forms a basis of $\Z^2$. Thus two bypass-related slopes form a lax basis. Conversely, for any $q/p, s/r$ forming a lax basis, $\Gamma_{q/p}, \Gamma_{s/r}$ are related by bypass surgery. So bypass-related slopes correspond precisely to lax bases. Moreover, a non-trivial bypass triple of sets of sutures of this type forms a lax superbasis; conversely, for any lax superbasis there exist sutures of these slopes forming a bypass triple.

This description of lax bases and superbases, or equivalently of bypass-related and bypass-triple slopes, forms the \emph{Farey graph}, which has vertices $\Q \cup \{\infty\}$ (i.e. slopes) and edges connecting lax bases. When drawn in the plane, with vertices $\Q \cup \infty$ lying on the circle $\R \cup \infty$, as a circle in the standard fashion, the Farey graph cuts the unit disc into triangles corresponding to lax superbases.

\begin{Proof}
We described above that for $m=-1,0,\infty$, a loop of slope $m$ in $H_1(T_1)$, and the suture element $c(\Gamma_m) \in V(T_1)^0$, have the same coordinates.

We first check that the same is also true for $m=1$. Sutures $\Gamma_{1/1}$ of slope $1$ are given by $\Omega_1 \Psi_{13} (xyxy)$ or $\Psi_1 \Omega_{31} (yxyx)$; from section \ref{sec:coherent_signs_for_gluings} these are both equal to $\Psi_1 \Omega_{11} (xy+yx)$, and hence have suture elements $\pm(1,1)$.

Thus, for $m=-1,0,1,\infty$, a loop of slope $m$ in $H_1(T_1)$, and the suture element $c(\Gamma_m) \in V(T_1)^0$, have the same coordinates.

The same argument applies any time an octagon homeomorphic to $O_{11}$ is glued to form $T_1$; the gluings need not be along the directions $(0,1)$ and $(1,0)$. Hence:
\begin{itemize}
\item
A bypass-related triple of sutures, with slopes forming a lax superbasis $\{\pm u, \pm v, \pm w\}$ of $H_1(T_1)$, has a triple of suture elements $\{c(\Gamma_u), c(\Gamma_v), c(\Gamma_w)\}$ forming a lax superbasis of $V(T_1)^0$.
\item
A bypass-related pair of sutures, with slopes forming a lax basis $\{\pm v, \pm w\}$ of $H_1(T_1)$, has a pair of suture elements $\{c(\Gamma_v), c(\Gamma_w)\}$ forming a lax basis of $V(T_1)^0$.
\item
Given a lax basis $\pm v, \pm w$ of $H_1(T_1)$, it includes into a lax superbasis in precisely two distinct ways $\{\pm v, \pm w, \pm (v+w)\}$ and $\{\pm v, \pm w, \pm(v-w)\}$; and correspondingly, the bypass-related pair of sutures $\Gamma_v, \Gamma_w$ includes into a bypass triple in precisely two distinct ways $\{\Gamma_v, \Gamma_w, \Gamma_{v+w}\}$ and $\{\Gamma_v, \Gamma_w, \Gamma_{v-w}\}$. The two corresponding triples of suture elements form the two distinct lax superbases containing the lax basis $\{c(\Gamma_v), c(\Gamma_w)\}$.
\end{itemize}
Thus, when the slopes $\pm u, \pm v, \pm (u+v) \in H_1(T_1)$, forming a lax superbasis, have the same coordinates as $c(\Gamma_u), c(\Gamma_v), c(\Gamma_w) \in V(T_1)^0$, the same is true for the lax superbasis $\pm u, \pm v, \pm(u-v)$. In other words, when the slopes of one superbasis containing $\pm u, \pm v$ have the same coordinates as their corresponding suture elements, the same is true for the other superbasis containing $\pm u, \pm v$. Now any primitive lax superbasis of $H_1(T_1)$ can be reached from $\pm(1,0), \pm(0,1), \pm (1,-1)$ by repeatedly restricting to a basis and extending to a superbasis; the dual of the Farey graph is connected.

Thus the set of sutures with an arc and loop of slope $m = \pm(p,q) \in H_1(T_1)$ has suture element $\pm(p,q) \in V(T_1)^0$ as desired.
\end{Proof}

\section{Sutured Floer homology}
\label{sec:SFH}

\subsection{Proofs of $SFH$ results}

We showed in \cite{Me10_Sutured_TQFT} that the following associations define a sutured TQFT satisfying axioms \ref{ax:1}--\ref{ax:10}:
\begin{itemize}
\item
Let $V(\Sigma,F) = SFH(- \Sigma \times S^1, - F \times S^1)$ with $\Z$ coefficients.
\item
A set of sutures $\Gamma$ on $(\Sigma, F)$ corresponds precisely to an isotopy class of contact structures $\xi$ on $\Sigma \times S^1$, such that the boundary $\partial \Sigma \times S^1$ is convex with dividing set $F \times S^1$ and positive/negative regions determined by the decomposition of $\partial \Sigma \backslash F$ into  positive and negative arcs $C_+ \cup C_-$ \cite{Gi00, GiBundles, Hon00II}. Let $c(\Gamma)$ be the \emph{contact invariant} $c(\xi) \subset SFH(-\Sigma \times S^1, -F \times S^1) = V(\Sigma,F)$ \cite{OSContact, HKMContClass, HKM09}. 
\item
For a gluing $\tau$ of the sutured background surface $(\Sigma,F)$, let $\Phi_\tau: SFH(- \Sigma \times S^1, - F \times S^1) \To SFH( - (\#_\tau \Sigma) \times S^1, -(\#_\tau F) \times S^1)$ be the map defined in \cite{HKM08} by the obvious inclusion of $\Sigma \times S^1 \hookrightarrow \#_\tau \Sigma \times S^1$, together with the canonical contact structure on $\#_\tau \Sigma \times S^1 - \Sigma \times S^1$ as convex neighbourhood of the boundary. In fact we can choose a sign on each $\Phi_\tau$ on each Euler class summand and take all $\Phi_\tau^i$ obtained by all possible choices of signs.
\end{itemize}

Now all of our results for sutured TQFT can immediately be applied to sutured TQFT.

\begin{Proof}[of theorem \ref{thm:isolating_vanishing}]
From the above, it is sufficient to prove that in sutured TQFT, for a set of sutures $\Gamma$ on a sutured background $(\Sigma, F)$, the following are equivalent:
\begin{enumerate}
\item
$c(\Gamma) \neq 0$.
\item
$c(\Gamma)$ is primitive.
\item
$\Gamma$ is not isolating.
\end{enumerate}
Proposition \ref{prop:nonzero} gives $(iii) \Rightarrow (ii)$. That $(ii) \Rightarrow (i)$ is obvious. And $(i) \Rightarrow (iii)$ (or rather, its contrapositive) is theorem \ref{thm:STQFT_isolating_vanishing}.
\end{Proof}

\begin{Proof}[of theorem \ref{thm:annulus_contact_elt_classification}]
Now immediate from the above and theorem \ref{thm:annulus_suture_elt_classification}.
\end{Proof}

\begin{Proof}[of theorem \ref{thm:punctured_torus_contact_elt_classification}]
Now immediate from the above and theorem \ref{thm:punctured_torus_suture_elt_classification}.
\end{Proof}

\subsection{Contact torsion and suture torsion}

Finally, in this section we prove theorem \ref{thm:torsion_vanishing}. First, however, we make some comments. The notions of torsion in contact topology, and torsion sutures in sutured TQFT, are similar but not identical.

In the above we have proved that the set of sutures $\Gamma_0$ on the annulus $(A, F_A)$ shown in figure \ref{fig:13} has suture element $0$. This corresponds to an $S^1$-invariant contact structure $\xi_0$ on $(A \times S^1, F_A \times S^1)$ which we now know has contact element $0$. But $\xi_0$ ``almost has $3\pi$-torsion'', which seems to be rather more than we need! On the other hand, the set of sutures on $(A,F_A)$ described by $\Phi_1(yx)$, which we shall now call $\Gamma_{2\pi}$, corresponding to a contact structure $\xi_{2\pi}$, ``almost has $2\pi$-torsion'', and yet has nonzero suture element.

To see why these manifolds ``almost'' have the torsion described, take a standard contact structure on $T^2 \times [0,1]$ with coordinates $((x,y),t)$ (which is also $A \times S^1$) with $3\pi$-torsion, i.e. $\xi = \ker ( \cos (3\pi t) \; dx - \sin (3\pi t) \; dy  )$. The boundary tori $T^2 \times \{0,1\}$ are then non-convex: they are pre-Lagrangian, with vertical characteristic foliation (i.e. in the $y$ direction). We can $C^\infty$ perturb these tori to make them convex; they can be taken then to have $2$ dividing curves, both vertical, so that the boundary dividing set is $F_A \times S^1$ as we desire. It's easy to check, then, that we have the contact structure $\xi_0$. Similarly, if we start with a standard $T^2 \times S^1$ with $2\pi$-torsion and boundary pre-Lagrangian tori with vertical characteristic foliation, we may again perturb the boundary tori to be convex with $2$ vertical dividing curves each; then we have the contact structure $\xi_{2\pi}$.

The contact structure $\xi_{2\pi}$ on $(A \times S^1, F_A \times S^1)$, with dividing set $\Gamma_{2\pi}$, is therefore $C^\infty$ close to a standard $2\pi$-torsion contact manifold, and yet has contact element nonzero. In perturbing the pre-Lagrangian torus boundary, we lose part of the torsion manifold. If, however, this manifold is enlarged in any way, for instance by a bypass attachment, then the resulting contact manifold with convex boundary has $2\pi$-torsion. In the context of sutured TQFT, we only consider contact manifolds with vertical dividing set on the boundary, and so to see a torsion contact structure we must go (``almost'') to $3\pi$-torsion. The simplest set of sutures in sutured TQFT corresponding to a contact structure with torsion is $\Gamma_0$ on $(A,F_A)$.

Thus, every set of sutures with torsion corresponds to a contact structure with torsion. And every contact structure with torsion on a contact manifold $\Sigma \times S^1$ with convex boundary and boundary dividing set of the form $F \times S^1$ corresponds to a set of sutures with torsion. So theorem \ref{thm:torsion_zero} is as good as result for torsion as we can imagine, within the ``contact geometry free'' subject of sutured TQFT. But not every contact structure with torsion is seen precisely in sutured TQFT. To prove theorem \ref{thm:torsion_vanishing}, then, we need a little contact geometry. 

As discussed at length in \cite{Hon00I}, tight contact toric annuli $T^2 \times I$ become ``larger'' as their boundary slopes change, moving around the ``circle'' of possible slopes ($\Q \cup \{\infty\}$ or $\R \cup \{\infty\}$). Torsion appears when the slopes of layers $T^2 \times \{\cdot\}$ have fully traversed a circle. If the boundaries are pre-Lagrangian, then we simply observe the slope of the characteristic foliation. If the boundaries are convex, we require that each boundary torus have a dividing set consisting of $2$ curves, and we then observe the slope of the dividing curves. As mentioned above (and as described explicitly in \cite[lemma 3.4]{Ghi05}), a pre-Lagrangian torus with slope $s$ can be perturbed to a convex torus with dividing set of slope $s$. Any enlargement of a contact manifold with convex boundary can be done via addition of bypasses. For a contact $T^2 \times [0,1]$ with convex boundary and boundary dividing set slope $s_1$ on $T^2 \times \{1\}$, adding a single bypass along $T^2 \times \{1\}$ gives a new boundary slope $s_2$ which forms a lax basis of $\Z^2$ with $s_1$; equivalently, $s_1, s_2 \in \Q \cup \{\infty\}$ are joined by an edge of the Farey graph. Such a bypass attachment can also be regarded as gluing onto the existing manifold a particular contact $T^2 \times I$ called a \emph{basic slice}. A basic slice is in this sense the simplest nontrivial tight contact $T^2 \times I$. We refer to \cite{Hon00I} for further details.

In any case, adding any bypass to $\xi_{2\pi}$ will produce a torsion contact structure; and any contact manifold nontrivially larger than $(T^2 \times I, \xi_{2\pi})$ contains a basic slice attached to $\xi_{2\pi}$. This immediately gives the following fact, which was observed by Massot in \cite{Massot09}.
\begin{lem}
\label{lem:torsion_basic_slice}
A contact structure $(M, \xi)$ has torsion if and only if it contains a submanifold contactomorphic to $\xi_{2\pi}$ with a basic slice attached.
\qed
\end{lem}

Denote the sutured toric annulus $T^2 \times [0,1]$ with two sutures on each boundary torus, of slopes $s_0$ on $T^2 \times \{0\}$ and $s_1$ on $T^2 \times \{1\}$, by $(T^2 \times I, s_0 \cup s_1)$, so $(A \times S^1, F_A \times S^1) = (T^2 \times I, \infty \cup \infty)$. As sutured manifolds, all basic slices are homeomorphic; a basic slice has two (isotopy classes of) tight contact structures, arising from adding a positive or negative bypass. As any tight contact structure $\xi$ on $(A \times S^1, F_A \times S^1)$ can be considered $S^1$-invariant, there is a contactomorphism of $\xi$ which sends any possible arc of attachment for a positive (resp. negative) bypass, to any other. Hence for a given contact structure $\xi$ on $(T^2 \times I, \infty \cup \infty)$, the addition of any positive (resp. negative) basic slice to $\xi$ produces a contactomorphic result.

Without loss of generality then, we consider attaching a basic slice to $(T^2 \times I, \infty \cup \infty)$ along $T^2 \times \{1\}$ to obtain slope $1$, i.e. we attach a $(T^2 \times I, \infty \cup 1)$ and obtain a $(T^2 \times I, \infty \cup 1)$. We will choose our orientations to be such that, adding bypasses or basic slices along $T^2 \times \{1\}$, slope \emph{increases}, so our basic slice $(T^2 \times I, \infty \cup 1)$ contains convex slices isotopic to $T^2 \times \{\cdot\}$ which have dividing set slope any given negative rational, or $0$, or positive rational less than $1$.

\begin{figure}
\centering
\includegraphics[scale=0.5]{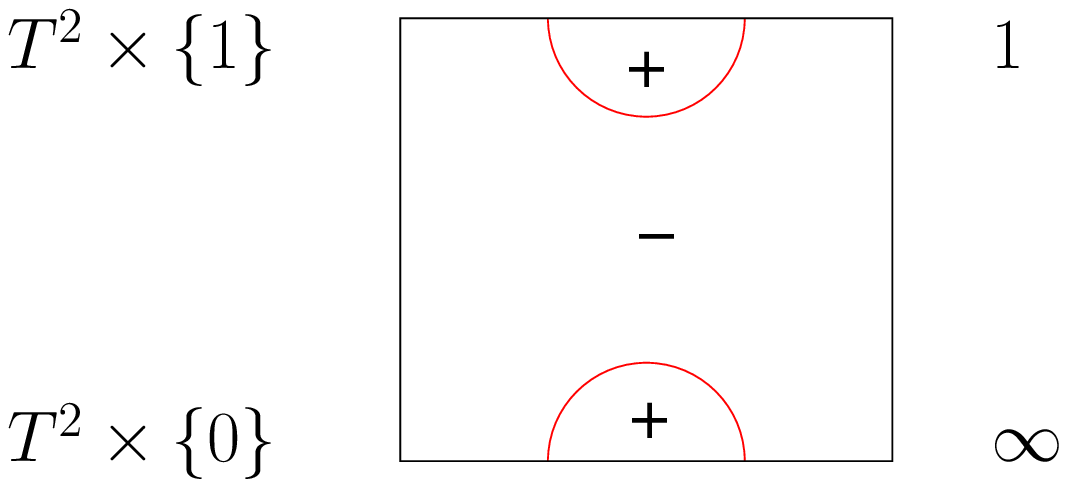}
\caption{Basic slice $(T^2 \times I, \infty \cup 1)$. The surface shown is the annulus $A = \{y=0\} \times I$; left and right sides are glued.} \label{fig:28}
\end{figure}

For definiteness we will consider this basic slice to be the attachment of a positive bypass. On the horizontal annulus $A = \{y=0\} \times I$, then, the dividing set is as shown in figure \ref{fig:28}. (One way to see why: consider attaching bypasses from either side to get to an intermediate torus with dividing set slope $0$.)

\begin{figure}
\centering
\includegraphics[scale=0.5]{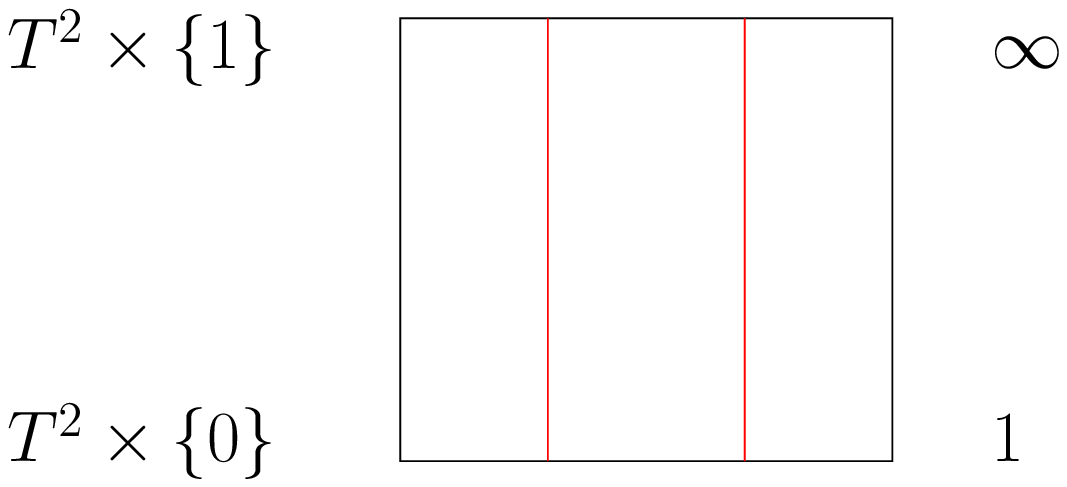}
\caption{Basic slice $(T^2 \times I, 1 \cup \infty)$. Again we show $A = \{y=0\} \times I$; left and right sides are glued.} \label{fig:29}
\end{figure}

We will also consider attaching a further basic slice $(T^2 \times I, 1 \cup \infty)$ onto $(T^2 \times I, \infty \cup 1)$ to obtain the original sutured manifold $(T^2 \times I, \infty \cup \infty)$.  Whether we take a positive or negative bypass, the dividing set on $A$ consists of two arcs, both running from one boundary of $A$ to the other, as in figure \ref{fig:29}. (One way to see why: for any other dividing set, we could find an intermediate torus isotopic to $T^2 \times \{ \cdot \}$ on which the slope must be $0$; a contradiction, since along this basic slice slopes increase from $1$ to $\infty$.)

\begin{figure}
\centering
\includegraphics[scale=0.5]{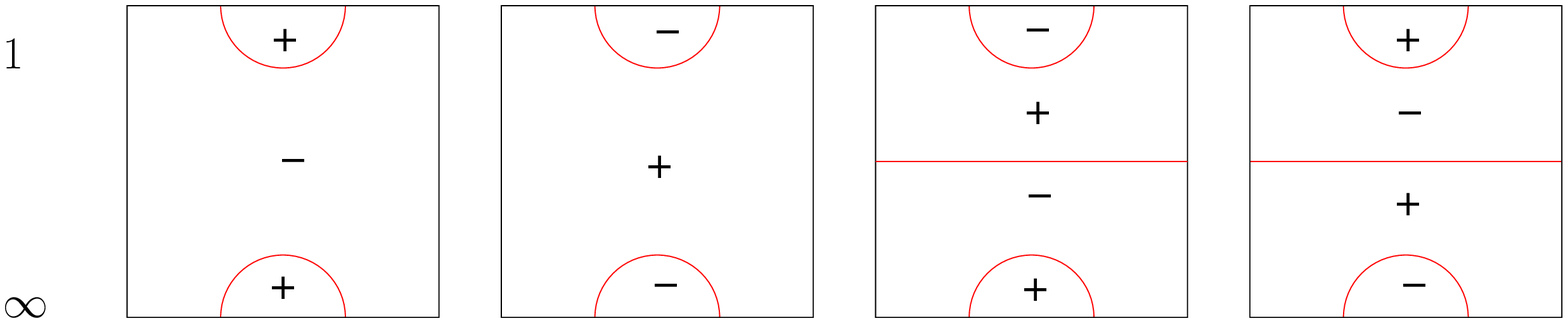}
\caption{Dividing sets on $A$ of the four contact structures on the basic slice $(T^2 \times I, \infty \cup 1)$ whose contact elements generate $SFH$. Again left and right sides are glued.} \label{fig:30}
\end{figure}

The sutured Floer homology of a basic slice was considered in \cite{HKM09}: Honda--Kazez--Mati\'{c} found there that (with $\Z$ coefficients) $SFH(T^2 \times I,\infty \cup 1) \cong \Z^4$, generated by the contact elements for four specific contact structures: the two basic slices (bypasses of either sign), and the two contact structures obtained by adding $\pi$ torsion to these. The four contact structures give dividing sets on $A = \{y=0\} \times I$ as shown in figure \ref{fig:30}. In particular, only one summand corresponds to contact structures with euler class evaluating to $2$ on $A$.

\begin{lem}
\label{lem:basic_slice_zero}
The contact element of the contact structure obtained from $\xi_{2\pi}$ by attaching a basic slice is $0$.
\end{lem}

We prove this from sutured TQFT by sandwiching the basic slice, which contains un-sutured-TQFT-like non-vertical sutures, between two manifolds which do have vertical sutures.
\begin{Proof}
We consider the two maps $J_1, J_2$ on $SFH$ obtained by the successive addition of basic slices described above, in order to add (``almost'') $\pi$ torsion to a contact structure.
\[
SFH(- T^2 \times I, - \infty \cup \infty) \stackrel{J_1}{\To} SFH(- T^2 \times I, - \infty \cup 1) \stackrel{J_2}{\To} SFH(-T^2 \times I, - \infty \cup \infty).
\]
Given a contact structure $\xi$ on any of these sutured manifolds, let $e_A(\xi)$ denote the euler class of $\xi$, evaluated on $A = \{y=0\} \times I$. Note that if $\xi$ is a contact structure on $(T^2 \times I, \infty \cup \infty)$ with $e_A = 0$, then after adding a positive basic slice as described above, we obtain a contact structure on $(T^2 \times I, \infty \cup 1)$ with $e_A = 2$ (this is clear in light of figure \ref{fig:28}). Then, adding another another basic slice, we obtain another contact structure with $e_A=2$ (in light of figure \ref{fig:29}). The relevant summands of sutured Floer homology are
\begin{align*}
\Z^2 &= SFH(-T^2 \times I, - \infty \cup \infty, e_A = 0) = V(A,F_A)^0,\\
\Z &= SFH(-T^2 \times I, -\infty \cup 1, e_A = 2),\\
\Z &= SFH(-T^2 \times I, - \infty \cup \infty, e_A = 2) = V(A,F_A)^2.
\end{align*}
This is clear from the above discussion and our earlier computations $V(A,F_A)^0 = \Z^2$, $V(A,F_A)^2 = \Z$. Restricting to these summands, we have
\[
\Z^2 \stackrel{J_1}{\To} \Z \stackrel{J_2}{\To} \Z.
\]
Now, as discussed above, $SFH(-T^2 \times I, -\infty \cup 1, e_A = 2)$ is generated by the contact element of the unique contact structure with $e_A=2$ in figure \ref{fig:30} (i.e. the leftmost), which is a basic slice; under $J_2$ this maps to the contact element of the ``almost $\pi$-torsion'' contact structure on $(A \times S^1, F \times S^1)$ obtained from joining the two basic slices; using figure \ref{fig:29}, the corresponding dividing set on $A$ consists of two boundary-parallel arcs enclosing positive discs. By theorem \ref{thm:annulus_contact_elt_classification}, this is a generator of the $e_A=2$ summand of $SFH(-T^2 \times I, - \infty \cup \infty)$. Thus $J_2$ gives an isomorphism $\Z \To \Z$.

On the other hand, in view of figures \ref{fig:28} and \ref{fig:29}, $J_1$ takes the contact element of the ``almost $2\pi$-torsion'' contact structure $\xi_{2\pi}$ corresponding to sutures $\Gamma_{2\pi}$ (which has $e_A = 0$), to the contact element of $\xi_{2\pi}$ plus a basic slice (which has $e_A = 2$); and then $J_2$ takes this to the contact element of the ``almost $3\pi$-torsion'' contact structure $\xi_0$ (which has $e_A = 2)$. Thus we have
\[
0 \neq c(\xi_{2\pi}) \stackrel{J_1}{\mapsto} c(\xi_{2\pi} \cup \text{basic slice}) \stackrel{J_2}{\mapsto} c(\xi_0) = 0.
\]
As discussed above, $c(\xi_{2\pi}) \neq 0$ but $c(\xi_0) = 0$. As $J_2$ is an isomorphism on the relevant summand, the contact element of $\xi_{2\pi}$ plus a basic slice is $0$.

The same result is obtained whether we add a positive or negative bypass, and add to either side of $T^2 \times I$.
\end{Proof}

\begin{Proof}[of theorem \ref{thm:torsion_vanishing}]
Suppose we have a contact manifold $(M, \xi)$ with a torsion contact structure. By lemma \ref{lem:torsion_basic_slice}, there is an embedded $T^2 \times I \subset M$ contactomorphic to $\xi_{2\pi}$ with a basic slice attached. By lemma \ref{lem:basic_slice_zero}, the contact structure on this $T^2 \times I$ has contact element $0$. From the inclusion $T^2 \times I \hookrightarrow M$, we obtain a map on $SFH$ which gives $c(\xi) = 0$.
\end{Proof}

The proof of this result by Ghiggini--Honda--Van Horn Morris shows that the contact invariant of a torsion toric annulus is zero; this is proved by considering the effect of Legendrian surgery on $SFH$ and the contact invariant. The proof of Massot (and also his proof that isolating dividing sets on $(\Sigma \times S^1, F \times S^1)$ give zero contact elements) uses the bypass relation and considers the effect of making various contact-geometric gluings to contact toric annuli. Our proof only uses two basic slices and their $SFH$ in order to sandwich a torsion $T^2 \times I$ between sutured TQFT constructions; from which inclusion maps give the result.

\addcontentsline{toc}{section}{References}

\small

\bibliography{danbib}
\bibliographystyle{amsplain}

\end{document}